\documentclass[12pt]{article}
\usepackage{amsfonts,amscd,a4}
\usepackage{color}
\usepackage{graphicx}
\usepackage{theorem}
\newtheorem{theo}{Theorem}[section]
\newtheorem{prop}[theo]{Proposition}
\newtheorem{lemma}[theo]{Lemma}
\newtheorem{coro}[theo]{Corollary}

{\theorembodyfont{\rm}
\newtheorem{example}[theo]{Example}
\newtheorem{remark}[theo]{Remark}}
\newcommand{\bA}{{\bf A}}

\newcommand{\cA}{{\mathcal A}}
\newcommand{\cB}{{\mathcal B}}
\newcommand{\cC}{{\mathcal C}}

\newcommand{\cF}{{\mathcal F}}
\newcommand{\cG}{{\mathcal G}}

\newcommand{\cI}{{\mathcal I}}
\newcommand{\cJ}{{\mathcal J}}
\newcommand{\cK}{{\mathcal K}}
\newcommand{\cL}{{\mathcal L}}
\newcommand{\cM}{{\mathcal M}}

\newcommand{\cO}{{\mathcal O}}

\newcommand{\cR}{{\mathcal R}}
\newcommand{\cS}{{\mathcal S}}
\newcommand{\cT}{{\mathcal T}}

\newcommand{\sC}{{\mathbb C}}

\newcommand{\sN}{{\mathbb N}}

\newcommand{\sT}{{\mathbb T}}

\newcommand{\sZ}{{\mathbb Z}}
\newcommand{\qed}{\rule{1ex}{1ex}}
\newcommand{\alg}{\mbox{\rm alg} \,}

\newcommand{\clos}{\mbox{\rm clos} \,}

\newcommand{\diag}{\mbox{\rm diag} \,}

\newcommand{\im}{\mbox{\rm im} \,}

\newcommand{\rank}{\mbox{\rm rank} \,}

\begin{document}
\title{Spatial discretization of Cuntz algebras}
\author{Steffen Roch}
\date{}
\maketitle
\begin{abstract}
The (abstract) Cuntz algebra $\cO_N$ is generated by non-unitary 
isometries and has therefore no intrinsic finiteness properties. 
To approximate the elements of the Cuntz algebra by finite-dimensional 
objects, we thus consider a spatial discretization of $\cO_N$ by 
the finite sections method. For we represent the Cuntz algebra as a 
(concrete) algebra of operators on $l^2(\sZ^+)$ and associate with 
each operator $A$ in this algebra the sequence $(P_n A P_n)$ of 
its finite sections. The goal of this paper is to examine the structure 
of the $C^*$-algebra $\cS(\cO_N)$ which is generated by all sequences 
of this form. Our main results are the fractality of a suitable 
restriction of the algebra $\cS(\cO_N)$ and a necessary and sufficient 
criterion for the stability of sequences in the restricted algebra. 
These results are employed to study spectral and pseudospectral 
approximations of elements of $\cO_N$. 
\end{abstract}
\section{Introduction} \label{s1}
Several classes of $C^*$-algebras are distinguished by intrinsic
finiteness properties. These properties can be used in principle
to approximate the elements of the algebra by finite-dimensional
(or discrete) objects and, thus, to discretize the algebra. Good
candidates for a discretization in that sense are AF-algebras and
quasidiagonal algebras. N. Brown \cite{Bro1} has pointed out that the discretization procedure works particularly well for irrational rotation algebras in which case the discrete approximations can not only be constructed effectively but also own excellent convergence properties (for example, the sequence of the approximations is fractal in a sense which will be explained
below).

At the other end of the scale there are $C^*$-algebras of infinite type which resist any intrinsic discretization. This fact
justifies to consider another kind of approximation of the
elements of a Cuntz algebra by finite rank operators, which we
call {\em spatial approximation} and which is based on the finite
sections method. Spatial approximation requires to represent the
algebra $^*$-isomorphically as a concrete $C^*$-algebra of
operators on an infinite dimensional Hilbert space $H$. All
concrete algebras considered in this paper will be separable; so
we can assume that $H$ is separable and, thus, isomorphic to the
Hilbert space $l^2(\sZ^+)$ of all sequences $(x_n)_{n \ge 0}$ of
complex numbers such that
\[
\|(x_n)\|^2 := \sum_{n \ge 0} |x_n|^2 < \infty.
\]
The sequences 
\[
e_i := (0, \, \ldots, \, 0, \, 1, \, 0, \, 0, \, \ldots)
\]
with the 1 standing at the $i$th position form a basis of $l^2(\sZ^+)$ to which we refer as the standard basis. For $n$ a positive integer, let $P_n$ denote the operator
\[
P_n : l^2(\sZ^+) \to l^2(\sZ^+), \quad (x_k)_{k \ge 0} \mapsto
(x_0, \, x_1, \, \ldots, \, x_{n-1}, \, 0, \, 0, \, \ldots),
\]
and set $P_0 := 0$. Thus, $P_n$ is the orthogonal projection of $l^2(\sZ^+)$ onto the span of the first $n$ elements of the standard basis.

For the finite sections method for a bounded linear operator $A$ on $l^2(\sZ^+)$, one replaces the equation $Au=f$ on $l^2(\sZ^+)$ by the sequence of the equations 
\begin{equation} \label{e91.1}
P_n A P_n u_n = P_n f, \quad n = 1, \, 2, \, 3, \, \ldots
\end{equation}
the solutions $u_n$ of which are sought in $\im P_n$. The sequence $(P_n A P_n)$ is called the sequence of the finite sections of $A$. This sequence is said to be {\em stable} if there is an $n_0$ such that the operators $P_n A P_n : \im P_n \to \im P_n$ are invertible for $n \ge n_0$ and if their inverses are uniformly bounded. If the sequence $(P_n A P_n)$ is stable, then the operator $A$ is invertible, the equations (\ref{e91.1}) possess unique solutions $u_n$ for all $n \ge n_0$ and for all right hand sides $f \in l^2(\sZ^+)$, and these solutions converge to the solution $u$ of the equation $Au=f$.

Let $\cA$ be a (for a moment not necessarily separable)
$C^*$-subalgebra of the algebra $L(l^2(\sZ^+))$ of all bounded
linear operators on $l^2(\sZ^+)$. We associate with each operator
$A$ in $\cA$ the sequence $(P_n A P_n)$ of its finite sections and consider the $C^*$-algebra $\cS(\cA)$ generated by this sequence.
Spatial discretization of the concrete algebra $\cA$ means to
study the associated algebra $\cS(\cA)$ of sequences. (To be
precise: we are only interested in asymptotic properties of the
sequences in $\cS(\cA)$ which are encoded in the quotient algebra
of $\cS(\cA)$ modulo sequences tending to zero in the norm. It
is this quotient algebra we are really interested in.)

Algebras of infinite type typically contain non-unitary
isometries. The perhaps simplest example, the universal algebra
generated by one isometry, is $^*$-isomorphic to the smallest
closed $^*$-subalgebra $\cT(C)$ of $L(l^2(\sZ^+))$ which contains
the operator
\[
V : l^2(\sZ^+) \to l^2(\sZ^+), \quad (x_k)_{k \ge 0} \mapsto (0,
\,  x_0, \, x_1, \, \ldots)
\]
of forward shift. This is the contents of a theorem by Coburn
\cite{Cob1}. The algebra $\cT(C)$ is also known as the {\em
Toeplitz algebra}, since each of its elements is of the form $T(c) + K$ where $T(c)$ is a Toeplitz operator with continuous
generating function $c$ and $K$ is a compact operator. The
structure of the associated algebra $\cS(\cT(C))$ (factored by the zero sequences) and, thus, the stability of the finite sections
method for operators in $\cT(C)$ are fairly well understood. For a detailed account on the finite sections method for Toeplitz
operators, presented in an algebraic language, see \cite{BSi2} and Section 1.4 in \cite{HRS2}.

The aim of the present paper is to go one step further and to examine the spatial discretization of algebras which are generated by a
finite number of non-commuting non-unitary isometries, namely the
Cuntz algebras. Recall that an {\em isometry} is an element $s$ of a unital $^*$-algebra for which $s^*s$ is the identity element.

Let $N \ge 2$. The {\em Cuntz algebra} $\cO_N$ is the universal
$C^*$-algebra generated by $N$ isometries $S_0, \, \ldots, \,
S_{N-1}$ with the property that
\begin{equation} \label{e91.2}
S_0 S_0^* + \ldots + S_{N-1} S_{N-1}^* = I.
\end{equation}
Cuntz algebras cannot be obtained as inductive limits of type I
$C^*$-algebras. In particular, they cannot be approximated by
finite dimensional algebras in the sense of $AF$-algebras. (For
these and other facts, consult Cuntz' pioneering paper
\cite{Cun1}. A nice introduction is also in \cite{Dav1}.) The
importance of Cuntz algebras in theory and applications cannot be
overestimated. Let me only mention Kirchberg's deep result that a
separable $C^*$-algebra is exact if and only if it embeds in the
Cuntz algebra $\cO_2$, and the role that representations of Cuntz
algebras play in wavelet theory and signal processing (see
\cite{Bla2,BrJ1} and the references therein).

To discretize the Cuntz algebra $\cO_n$ by the finite sections
method, we have to represent this algebra as a $C^*$-subalgebra of
$L(l^2(\sZ^+))$. Since Cuntz algebras are simple, every
$C^*$-algebra which is generated by $N$ isometries $S_0, \,
\ldots, \, S_{N-1}$ which fulfill (\ref{e91.2}) is $^*$-isomorphic
to $\cO_N$. Thus, $\cO_N$ is $^*$-isomorphically to the smallest
$C^*$-subalgebra of $L(l^2(\sZ^+))$ which contains the operators
\begin{equation} \label{e91.3}
S_i : (x_k)_{k \ge 0} \mapsto (y_k)_{k \ge 0}
\quad \mbox{with} \quad
y_k := \left\{
\begin{array}{ll}
x_r & \mbox{if} \; k = rN + i \\
0   & \mbox{else}
\end{array}
\right.
\end{equation}
for $i = 0, \, \ldots, \, N-1$. We denote the (concrete) Cuntz
algebra generated by the operators $S_i$ in (\ref{e91.3}) also by
$\cO_N$. It is this concrete Cuntz algebra for which we will
examine the sequence algebra $\cS(\cO_N)$ (modulo zero sequences)
in what follows. One should mention at this point that the
abstract Cuntz algebra $\cO_N$ has an uncountable set of
equivalence classes of irreducible representations.
Representations of $\cO_N$ different from (\ref{e91.3}) will
certainly lead to different algebras $\cS(\cO_N)$. The relations
between these algebras are not yet clearly understood.

The paper is organized as follows. In Section \ref{s2} we collect
some basic facts from numerical analysis, centered around the
notions of stability and fractality. In Section \ref{s3} we
examine the full algebra of the finite sections method for
$\cO_N$. We show that this algebra fails to be fractal. This
observation suggests to pass to a restricted algebra of finite
sections sequences, and it is in fact this restricted algebra
which will be studied in what follows. In Section \ref{s4} we
divide the stability problem into an invertibility problem in a
certain quotient algebra of $\cS(\cO_N)$ and a lifting problem for
an ideal of $\cS(\cO_N)$. The first of these problems is also
solved in Section \ref{s4}, whereas the solution of the lifting
problem is subject of Section \ref{s6}. The treatment of the
lifting problem is prepared by Section \ref{s5} where we consider
an algebra of {\em stratified Toeplitz operators} related to the
Cuntz algebra and examine its structure. As already mentioned, in
Section \ref{s6} we finish the proof of our stability result.
These results are then applied to examine spectral approximation
and Fredholm properties of sequences in (the restriction of)
$\cS(\cO_N)$. We will see that a sequence in this algebra is
Fredholm if and only if it is stable, which has remarkable
consequences for the asymptotic behavior of the singular values
of the $n$th approximant when $n$ tends to infinity.

Coburn's already mentioned result suggests to consider the
Toeplitz algebra as the Cuntz algebra $\cO_1$. But one should have in mind that the main properties of $\cO_1$ and of $\cO_N$ for $N
> 1$ are completely different from each other. For example, the compact operators $K(l^2(\sZ^+))$ form a closed ideal of $\cO_1$, and the quotient $\cO_1/K(l^2(\sZ^+))$ is isomorphic to $C(\sT)$, whereas $\cO_N$ is simple if $N \ge 2$. These differences continue to the corresponding sequence algebras $\cS(\cO_1)$ and $\cS(\cO_N)$ for $N > 1$. A main point is that $\cS(\cO_1)$ contains an ideal which is constituted of two exemplars of the ideal $K(l^2(\sZ^+))$, and the irreducible representations, $W_1$ and $W_2$ say, of $\cS(\cO_1)$ which come from this ideal are {\em sufficient} in the sense that a sequence $\bA = (A_n)$ in $\cS(\cO_1)$ is stable if and only if $W_1(\bA)$ and $W_2(\bA)$ are invertible. It turns out that this fact implies an effective criterion to check the stability of a sequence in $\cS(\cO_1)$. In contrast to this, if $N > 1$, then $\cS(\cO_N)$ has only one non-trivial ideal. We will construct an injective representation of this ideal, and will then observe that this representation extends to a representation, $W$ say, of $\cS(\cO_N)$ which is injective on all of $\cS(\cO_N)$. Thus, roughly speaking, our stability result will say that a sequence $\bA$ in $\cS(\cO_N)$ is stable if and only if the operator $W(\bA)$ is invertible. At the first glance, this result might seem to be useless since the stability of $\bA$ is not easier to check than the invertibility of $W(\bA)$. So why this effort, if many canonical homomorphisms on $\cS(\cO_N)$ own the same property as $W$: the identical mapping and the faithful representation via the GNS-construction, for example. What is important is the concrete form of the mapping $W$ constructed below: it is defined by means of strong limits of operator sequences, and this special form implies an immediate proof of the fractality of the (restricted) algebra $\cS(\cO_N)$.
\section{Preliminaries from numerical analysis} \label{s2}
\subsection{Algebras of matrix sequences} \label{ss21}
Let $\cF$ denote the set of all bounded sequences $\bA = (A_n)$ of matrices $A_n \in \sC^{n \times n}$. Equipped with the
operations
\[
(A_n) + (B_n) := (A_n + B_n), \quad
(A_n) (B_n) := (A_n B_n), \quad (A_n)^* := (A_n^*)
\]
and the norm 
\[
\|\bA\|_\cF := \|A_n\|,
\]
the set $\cF$ becomes a $C^*$-algebra, and the set $\cG$ of all sequences $(A_n) \in \cF$ with $\lim \|A_n\| = 0$ forms a closed ideal of $\cF$. The relevance of the algebra $\cF$ and its ideal $\cG$ in our context stems from the fact (following from a simple Neumann series argument) that a sequence $(A_n) \in \cF$ is stable if, and only if, the coset $(A_n) + \cG$ is invertible in the quotient algebra $\cF/\cG$. Thus, every stability problem is equivalent to an invertibility problem in a suitably chosen $C^*$-algebra.

Let further stand $\cF^C$ for the set of all sequences $\bA = (A_n)$ of operators $A_n : \im P_n \to \im P_n$ with the property that the sequences $(A_nP_n)$ and $(A_n^*P_n)$ converge strongly. By the uniform boundedness principle, the quantity $\sup \|A_n P_n\|$ is finite for every sequence $\bA$ in $\cF^C$. This, if we identify each operator $A_n$ on $\im P_n$ with its matrix representation with respect to the basis $e_0, \, \ldots, \, e_{n-1}$ of $\im P_n$, we can consider $\cF^C$ as a closed and symmetric subalgebra of $\cF$ which contains $\cG$ as its ideal. 
Note that the mapping
\begin{equation} \label{e91.5}
W : \cF^C \to L(l^2(\sZ^+)), \quad (A_n) \mapsto \mbox{s-lim} \,
A_nP_n
\end{equation}
is a $^*$-homomorphism.
\subsection{Discretization of concrete algebras} \label{ss22}
Let $\cA$ be a (not necessarily separable) $C^*$-subalgebra of the algebra $L(l^2(\sZ^+))$. We write $D$ for the mapping of spatial
(= finite sections) discretization, i.e.,
\begin{equation} \label{e91.6}
D : L(l^2(\sZ^+)) \to \cF, \quad A \mapsto (P_n A P_n),
\end{equation}
and we let $\cS(\cA)$ stand for the smallest closed
$C^*$-subalgebra of the algebra $\cF$ which contains all sequences $D(A)$ with $A \in \cA$. Clearly, $\cS(\cA)$ lies even in $\cF^C$, and the mapping $W$ in (\ref{e91.5}) induces a $^*$-homomorphism from $\cS(\cA)$ onto $\cA$. On this level, one cannot say much about algebra $\cS(\cA)$. The little one can say will follow easily from the following simple facts.
\begin{prop} \label{p91.7}
Let $\cA$ and $\cB$ be $C^*$-algebras, $D : \cA \to \cB$ a linear
contraction, and $W : \cB \to \cA$ a $C^*$-homomorphism such that
$W(D(A)) = A$ for every $A \in \cA$. Then \\[1mm]
$(a)$ $D$ is an isometry, $D(\cA)$ is a closed linear subspace of
$\cB$, and $\alg D(\cA)$, the smallest closed $C^*$-subalgebra of
$\cB$ which contains $D(\cA)$, splits into the direct sum
\begin{equation} \label{e91.8}
\alg D(\cA) = D(\cA) \oplus (\ker W \cap \alg D(\cA)).
\end{equation}
Moreover, for every $A \in \cA$,
\begin{equation} \label{e91.9}
\|D(A)\| = \min_{K \in \ker W} \|D(A) + K\|.
\end{equation}
$(b)$ If $\cB = \alg D(\cA)$, then $\ker W$ coincides with the
quasicommutator ideal of $\cB$, i.e., with the smallest closed
ideal of $\cB$ which contains all quasicommutators $D(A_1) D(A_2)
- D(A_1 A_2)$ with $A_1, \, A_2 \in \cA$.
\end{prop}
{\bf Proof.} $(a)$ Let $A \in \cA$. The inequality
\[
\|A\| = \|W(D(A))\| \le \|D(A)\| \le \|A\|
\]
shows that $D$ is an isometry; hence, $D(\cA)$ is a closed
subspace of $\cB$. Let $B \in D(\cA) \cap \ker W$. Write $B =
D(A)$ with $A \in \cA$. From $W(B) = 0$ we get $A = W(D(A)) = W(B)
= 0$, whence $B = 0$. Thus, $D(\cA) \cap \ker W = \{0\}$.

Let $B \in \alg D(\cA)$. Then $W(B - D(W(B))) = W(B) - W(B) = 0$,
hence
\[
B = D(W(B)) + (B - D(W(B)) \in D(\cA) + \ker W,
\]
whence $\alg D(\cA) = D(\cA) + (\ker W \cap \alg D(\cA))$. This
proves (\ref{e91.8}). To check (\ref{e91.9}), let $A \in \cA$ and
$K \in \ker W$. Then
\[
\|A\| = \|W(D(A) + K)\| \le \|D(A) + K\|
\]
which implies that $\|D(A)\| \le \|D(A) + K\|$ since $D$ is an
isometry. \\[2mm]
$(b)$ Since $W$ is a homomorphism and $W \circ D$ is the identity
on $\cA$, one has $D(A_1)D(A_2) - D(A_1A_2) \in \ker W$ for all
$A_1, \, A_2 \in \cA$. Thus, $\ker W$ contains the quasicommutator ideal. For the reverse inclusion, let $K \in \ker W$ and $n$ a
positive integer. Since $K \in \alg D(\cA)$, there are sums of
products
\[
K_n = \sum \prod D(A_{ij}^{(n)}) \quad \mbox{with} \quad
A_{ij}^{(n)} \in \cA
\]
such that $\|K - K_n\| \le 1/n$. Clearly, each $K_n$ can be
written as
\[
K_n = D(\sum \prod A_{ij}^{(n)}) + Q_n
\]
with an element $Q_n$ in the quasicommutator ideal. From
\[
\|K - Q_n\| \le \|K - K_n\| + \|D(\sum \prod A_{ij}^{(n)})\|
\]
and
\[
\|D(\sum \prod A_{ij}^{(n)})\| \le \|\sum \prod A_{ij}^{(n)}\| =
\|W(K_n)\| = \|W(K_n - K)\| \le \|K_n - K\|
\]
we conclude that $\|K - Q_n\| \le 2 \|K - K_n\| \le 1/n$. Thus,
$K$ can be approximated as closely as desired by elements in the
quasicommutator ideal. Since the quasicommutator ideal is closed,
the assertion follows. \hfill \qed \\[3mm]
We apply the preceding proposition in the following context:
\begin{itemize}
\item
$\cA$ is a $C^*$-subalgebra of $L(l^2(\sZ^+))$,
\item
$\cB = \cS(\cA)$,
\item
$D$ is the restriction of the discretization (\ref{e91.6}) to
$\cA$, and
\item
$W$ is the restriction of the homomorphism (\ref{e91.5}) to
$\cS(\cA)$.
\end{itemize}
Then Proposition \ref{p91.7} specializes to the following.
\begin{prop} \label{p91.10}
Let $\cA$ be a $C^*$-subalgebra of $L(l^2(\sZ^+))$. Then the
finite sections discretization $D : \cA \to \cF$ is an isometry,
and $D(\cA)$ is a closed subspace of $\cS(\cA)$. The algebra
$\cS(\cA)$ splits into the direct sum
\[
\cS(\cA) = D(\cA) \oplus (\ker W \cap \cS(\cA)),
\]
and one has
\[
\|D(A)\| = \min_{K \in \ker W} \|D(A) + K\|
\]
for every operator $A \in \cA$. Finally, $\ker W \cap \cS(\cA)$ is
equal to the quasicommutator ideal of $\cS(\cA)$, i.e., to the
smallest closed ideal of $\cS(\cA)$ which contains all sequences
$(P_n A_1 P_n A_2 P_n - P_n A_1 A_2 P_n)$ with operators $A_1, \,
A_2 \in \cA$.
\end{prop}
We denote the ideal $\ker W \cap \cS(\cA)$ of $\cS(\cA)$ by
$\cJ(\cA)$. Since the first item in the decomposition $D(\cA)
\oplus \cJ(\cA)$ of $\cS(\cA)$ is isomorphic (as a linear space)
to $\cA$, a main part of the description of the algebra $\cS(\cA)$
is to identify the ideal $\cJ(\cA)$.
\begin{remark}
Blackadar and Kirchberg gave an abstract characterization of
$C^*$-subalgebras of $\cF/\cG$ as generalized inductive limits of
finite-dimensional algebras, which they call MF-algebras (see
the last sections of Blackadar's monograph \cite{Bla2}). In this
sense, the ideal $\cJ(\cA)$ and its image $\cJ(\cA)/\cG$ in
$\cF/\cG$ can be considered as a measure for the deviation of an
algebra from being an MF-algebra. For example, for the Toeplitz
algebra $\cA = \cT(C)$ which is generated by one non-unitary
isometry, one knows that $\cJ(\cA)/\cG$ is $^*$-isomorphic to the
ideal $K(l^2(\sZ^+))$ of the compact operators.

Note that since $\cS(\cA)/\cJ(\cA)$ is canonically $^*$-isomorphic to $\cA$, the above construction implies a simple proof of the
(well-known) fact that every $C^*$-subalgebra of $L(l^2(\sZ^+))$
is $^*$-isomorphic to a quotient of an MF-algebra. \hfill \qed
\end{remark}
\subsection{Fractal algebras of matrix sequences}
\label{ss23}
The mapping which associates with every subalgebra $\cA$ of
$L(l^2(\sZ^+))$ the sequence algebra $\cS(\cA)$ does not mind the
individual properties of $\cA$. If $\cA$ is separable, this fact
can be compensated to some extent by passing to a {\em fractal
restriction of} $\cS(\cA)$, which is defined as follows.

Let $\eta : \sN \to \sN$ be a strongly monotonically increasing
sequence. By $\cF_\eta$ we mean the set of all subsequences
$(A_{\eta(n)})$ of sequences $(A_n)$ in $\cF$. As in Section
\ref{ss21}, $\cF_\eta$ can be made to a $C^*$-algebra in a natural
way. The $^*$-homomorphism
\[
R_\eta : \cF \to \cF_\eta, \quad (A_n) \mapsto (A_{\eta(n)})
\]
is called the restriction of $\cF$ onto $\cF_\eta$. It maps the
ideal $\cG$ of $\cF$ onto a closed ideal $\cG_\eta$ of $\cF_\eta$. For every subset $\cS$ of $\cF$, we abbreviate $R_\eta \cS$ by
$\cS_\eta$.

Let $\cS$ be a $C^*$-subalgebra of $\cF$. A $^*$-homomorphism $W$
from $\cS$ into a $C^*$-algebra $\cB$ is called {\em fractal} if,
for every strongly monotonically increasing sequence $\eta : \sN
\to \sN$, there is a mapping $W_\eta : \cS_\eta \to \cB$ such that $W = W_\eta R_\eta|_\cS$. A $C^*$-subalgebra $\cS$ of $\cF$ is
called {\em fractal}, if the canonical homomorphism
\[
\cS \to \cS/(\cS \cap \cG), \quad \bA \mapsto \bA + (\cS \cap \cG)
\]
is fractal.

Thus, if $\cS$ is a fractal algebra, then every sequence in $\cS$
can be uniquely rediscovered from each of its (infinite)
subsequences up to a sequence tending to zero. In that sense, the
essential information on a sequence in $\cS$ is already stored in
each of its subsequences. These algebras were called {\em fractal} in \cite{RoS5} in order to emphasize exactly this 
self-similarity aspect.

If $\cA$ is a $C^*$-subalgebra of $L(l^2(\sZ^+))$ and $\eta : \sN
\to \sN$ a strongly monotonically increasing sequence, then we
denote the restriction $R_\eta (\cS(\cA))$ of the associated
sequence algebra by $\cS_\eta(\cA)$. The following was proved in
\cite{Roc0}.
\begin{theo} \label{t91.11}
Let $\cA$ be a separable $C^*$-subalgebra of $L(l^2(\sZ^+))$. Then
there exists a strongly monotonically increasing sequence $\eta :
\sN \to \sN$ such that the restricted algebra $\cS_\eta(\cA)$ is
fractal.
\end{theo}
Since the sequence $\eta$ depends on $\cA$, the fractal
restrictions $\cS_\eta(\cA)$ of $\cS(\cA)$ will reflect the
structure of $\cA$ in a much higher extent than the full algebra
$\cS(\cA)$ of the finite sections discretization.

One should mention that there are further good reasons to be
interested in fractal algebras. One of them is that sequences in
fractal algebras exhibit an excellent asymptotic behavior. For
example, if $(A_n)$ is a sequence in a fractal algebra, then
several important spectral quantities of the $A_n$ (e.g., the set
of the singular values, the pseudospectrum, the numerical range)
converge with respect to the {\em Hausdorff metric} as $n$ tends to infinity (see Chapter 3 in \cite{HRS2} and Section \ref{ss64}
below). Another reason came up in \cite{Roc9} where we observed
that the property of fractality determines the ideal structure
of the algebra to a large extent.
\section{The full algebra $\cS(\cO_N)$} \label{s3}
\subsection{The full algebra of the finite sections method for
operators in $\cO_N$} \label{ss31}
In accordance with the above notations, let $\cS(\cO_N)$ denote
the smallest closed subalgebra of $\cF$ which contains all
sequences $(P_n A P_n)$ with $A$ in the concrete Cuntz algebra
$\cO_N$. Since $(P_n A P_n)^* = (P_n A^* P_n)$, $\cS(\cO_N)$ is a
$C^*$-algebra. The isometries $S_i$ are defined as in
(\ref{e91.3}). We further abbreviate $\Omega := \{0, \, 1, \,
\ldots, \, N-1\}$.
\begin{lemma} \label{l92.1}
$\cS(\cO_N)$ is the smallest $C^*$-subalgebra of $\cF$ which
contains all sequences $(P_n S_j P_n)$ with $j \in \Omega$.
\end{lemma}
{\bf Proof.} For a moment, let $\cS^\prime$ denote the smallest
closed and symmetric subalgebra of $\cF$ which contains all
sequences $(P_n S_j P_n)$ with $j \in \Omega$. Evidently,
$\cS^\prime \subseteq \cS(\cO_N)$. For the reverse inclusion, note
first that
\begin{equation} \label{e92.2}
S_i^* S_j = 0 \quad \mbox{whenever} \; i \neq j.
\end{equation}
Indeed, this follows for the operators (\ref{e91.3}) by
straightforward calculation, but it also follows easily from the
Cuntz axiom (\ref{e91.2}): Multiply (\ref{e91.2}) from the left by
$S_i^*$ and from the right by $S_j$ and take into account that a
sum of positive elements in a $C^*$-algebra is zero if and only if
each of the elements is zero.

From (\ref{e92.2}) we conclude that every finite word with letters
in the alphabet $\{ S_1, \, \ldots, \, S_N, \, S_1^*, \, \ldots,
S_N^* \}$ is of the form
\begin{equation} \label{e92.3}
S_{i_1} S_{i_2} \ldots S_{i_k} S^*_{j_1} S^*_{j_2} \ldots S^*_{j_l}
\quad \mbox{with} \quad i_s, \, j_t \in \Omega
\end{equation}
(Lemma 1.3 in \cite{Cun1}). Further one easily checks that
\begin{equation} \label{e92.4}
P_n S_j = P_n S_j P_n \quad \mbox{and} \quad S_j^* P_n = P_n S_j^*
P_n
\end{equation}
for every $j \in \Omega$ and every $n \in \sN$. Thus, if $A$ is
any word of the form (\ref{e92.3}), then
\[
P_n A P_n = P_n S_{i_1} P_n \cdot P_n S_{i_2} P_n \ldots P_n S_{i_k} P_n
\cdot P_n S^*_{j_1} P_n \cdot P_n S^*_{j_2} P_n \ldots P_n S^*_{j_l} P_n
\in \cS^\prime.
\]
Since the set of all linear combinations of the words
(\ref{e92.3}) is dense in $\cO_N$, it follows that $\cS(\cO_N)
\subseteq \cS^\prime$.
\hfill \qed \\[3mm]
Recall that an element $S$ of a $C^*$-algebra is called a partial
isometry if $SS^*S = S$. If $S$ is a partial isometry, then $SS^*$
and $S^*S$ are projections (i.e., self-adjoint idempotents),
called the {\em range projection} and the {\em initial projection}
of $S$, respectively. Conversely, if $S^*S$ (or $SS^*$) is a
projection for an element $S$, then $S$ is a partial isometry.
Recall also that projections $P$ and $Q$ are called orthogonal if
$PQ = 0$.
\begin{lemma} \label{l92.5}
Every sequence $(P_n S_i P_n)$, $i \in \Omega$, is a partial
isometry in $\cF$, and the corresponding range projections are
orthogonal if $i \neq j$. Moreover,
\begin{equation} \label{e92.6}
P_n S_i^* P_n S_j P_n = 0 \quad \mbox{if} \; \; i \neq j,
\end{equation}
and
\begin{equation} \label{e92.7}
P_n S_0 P_n S_0^* P_n + \ldots + P_n S_{N-1} P_n S_{N-1}^* P_n = P_n.
\end{equation}
\end{lemma}
{\bf Proof.}  The identities (\ref{e92.4}) imply that
\begin{equation} \label{e92.8}
P_n S_i S_i^* P_n = P_n S_i P_n S_i^* P_n
\end{equation}
for every $i \in \Omega$ and every $n \in \sZ^+$. The operators
$S_i S_i^*$ are projections, and their matrices with respect to
the standard basis of $l^2(\sZ^+)$ are of diagonal form. Hence,
the projections $S_i S_i^*$ and $P_n$ commute, which implies that
the left-hand side of (\ref{e92.8}) is a projection. Hence, $(P_n
S_i P_n)$ is a partial isometry in $\cF$, and $(P_n S_i S_i^*
P_n)$ is the associated range projection.

Let $i \neq j$ be in $\Omega$. The fact that the $P_n$ and the
$S_i S_i^*$ commute further implies together with (\ref{e92.2})
that
\[
(P_n S_i S_i^* P_n) (P_n S_j S_j^* P_n) = P_n S_i S_i^* S_j S_j^* P_n
= 0.
\]
Multiplying $P_n S_i S_i^* P_n S_j S_j^* P_n = 0$ from the left by
$P_n S_i^* P_n$ and from the right by $P_n S_j P_n$ yields
(\ref{e92.6}). Finally, (\ref{e92.7}) follows by summing up the
equalities (\ref{e92.6}) over $i \in \Omega$ and from axiom
(\ref{e91.2}). \hfill \qed \\[3mm]
Thus, the generating sequences $(P_n S_i P_N)$, $i \in \Omega$,
are still subject of the Cuntz axiom (\ref{e91.2}), but note they
are partial isometries only and no longer isometries (which is not
a surprise since the algebra $\cF$, being a product of finite
dimensional algebras, cannot contain non-unitary isometries).

The first assertion of the preceding lemma holds more generally.
\begin{lemma} \label{l92.9}
Let $i = (i_1, \, i_2, \, \ldots, \, i_k) \in \Omega^k$. Every
product
\[
(P_n S_{i_1} P_n S_{i_2} P_n \ldots P_n S_{i_k} P_n)
\]
is a partial isometry in $\cF$.
\end{lemma}
Indeed, from (\ref{e92.4}) we conclude that
\[
P_n S_{i_1} \ldots S_{i_k} S_{i_k}^* \ldots S_{i_1}^* P_n = P_n
S_{i_1} P_n \ldots P_n S_{i_k} P_n P_n S_{i_k}^* P_n \ldots P_n
S_{i_1}^* P_n.
\]
Since $S_{i_1} \ldots S_{i_k} S_{i_k}^* \ldots S_{i_1}^*$ is a
projection of diagonal form, the assertion follows as in Lemma
\ref{l92.5}. \hfill \qed
\subsection{Initial projections} \label{ss32}
For $i = (i_1, \, i_2, \, \ldots, \, i_k) \in \Omega^k$,
abbreviate $S_i := S_{i_1} S_{i_2} \ldots S_{i_k}$. By Lemma
\ref{l92.9}, the sequence
\[
(P_n S_i P_n) = (P_n S_{i_1} P_n S_{i_2} P_n \ldots P_n S_{i_k} P_n)
\]
is a partial isometry. We are going to determine its initial
projection. The result will indicate that it is more natural to
consider a certain restriction of the finite sections algebra
$\cS(\cO_N)$ rather than the full algebra $\cS(\cO_N)$.

For every real number $x$, let $\{x\}$ denote the smallest integer
which is greater than or equal to $x$.
\begin{prop} \label{p92.10}
For $i = (i_1, \, i_2, \, \ldots, \, i_k) \in \Omega^k$, set
$v_{i, k} := i_1 + i_2 N + \ldots + i_k N^{k-1}$. Then
\begin{equation} \label{e92.11}
P_n S_i^* S_i P_n = P_{\{ (n - v_{i,k})/N^k \} }.
\end{equation}
\end{prop}
{\bf Proof.} It follows from the definition of $S_i$ that there
are numbers $v_{i,k}$ and $d_k$ such that
\begin{eqnarray*}
\lefteqn{S_i = S_{i_1} \ldots S_{i_k} : \quad (x_k)_{k \ge 0}
\mapsto} \\
&& ( \underbrace{0, \, \ldots, \, 0,}_{v_{i,k}} \, x_0, \,
\underbrace{0, \, \ldots, \, 0,}_{d_k} \, x_1, \,
\underbrace{0, \, \ldots, \, 0,}_{d_k} \, x_2, \, \ldots).
\end{eqnarray*}
The initial values
\[
v_{(i_1), 1} = i_1   \quad \mbox{and} \quad d_1 = N-1
\]
together with the recursions
\[
i_1 + v_{(i_2, \, \ldots, \, i_{k+1}), k} \cdot N =
v_{(i_1, \, \ldots, \, i_{k+1}), k+1} \quad \mbox{and} \quad
N-1 + d_k \cdot N = d_{k-1}
\]
imply via induction that
\begin{equation}
v_{i,k} = i_1 + i_2 N + \ldots + i_k N^{k-1} \quad \mbox{and}
\quad d_k = N^k-1.
\end{equation}
It follows that the $j$th component of $S_i^*$ applied to $x =
(x_k)_{k \ge 0}$ is
\[
(S_i^* x)_j = x_{v_{i,k} + j (d_k + 1)} = x_{v_{i,k} + j N^k}.
\]
Consequently,
\[
P_n S_i^* P_n S_i P_n = P_j \quad \mbox{if} \quad
v_{i,k} + (j-1) N^k < n \le v_{i,k} + j N^k,
\]
whence $j = \{ (n - v_{i,k})/N^k \}$. \hfill \qed
\subsection{The need of restrictions} \label{ss33}
We specialize the result of Proposition \ref{p92.10} to the case
$k = 1$ and consider two examples. If $n = jN$ is a multiple of
$N$, then the initial projections of $P_n S_i P_n$ are independent
of $i$. Indeed, from
\[
\{ (n-i)/N \} = \{ (j N-i)/N \} = \{ j - i/N \} = j
\]
we obtain
\begin{equation} \label{e92.12}
P_{jN} S_i^* P_{jN} S_i P_{jN} = P_j \quad \mbox{for all} \; i \in
\Omega.
\end{equation}
On the other hand, one has
\begin{equation} \label{e92.13}
P_n S_0^* P_n S_0 P_n - P_n S_1^* P_n S_1 P_n = \left\{
\begin{array}{ll}
P_{j+1} - P_j & \mbox{if} \; n = jN+1, \\
0             & \mbox{else}.
\end{array}
\right.
\end{equation}
Thus, the sequence
\begin{equation} \label{e92.14}
(P_n S_0^* P_n S_0 P_n - P_n S_1^* P_n S_1 P_n)_{n \ge 1}
\end{equation}
possesses both a subsequence consisting of zeros only (take
$\eta(n) := nN$) and a subsequence consisting of non-zero
projections (if $\eta(n) := nN + 1$). This shows that the algebra
$\cS(\cO_N)$ cannot be fractal.

Moreover, a closer look reveals that also the restricted algebra
$\cS_\eta(\cO_N)$, with $\eta(n) := nN$, cannot be fractal: For
the sequence
\begin{equation} \label{e92.15}
(P_{nN} (S_0^*)^2 P_{nN} S_0^2 P_{nN} -
P_{nN} (S_1^*)^2 P_{nN} S_1^2 P_{nN})_{n \ge 1}
\end{equation}
one observes a similar unpleasant behavior as before. We will see
later on that
\begin{equation} \label{e92.16}
\eta(n) := N^n
\end{equation}
is the correct choice for the restriction $\eta$, since it will
indeed guarantee the fractality of the restricted algebra
$\cS_\eta(\cO_N)$.

It is interesting to observe that there are at least two further
arguments which also suggest the choice (\ref{e92.16}) for the
restriction $\eta$. The first one comes from asymptotic numerical
analysis again. Consider the set of all sequences $(A_n)$ in $\cF$
with
\[
\sup_{n \ge 0} \rank A_n < \infty.
\]
The closure of this set in $\cF$ is a two-sided ideal $\cK$ of
$\cF$, the elements of which are called {\em compact sequences}.
The ideal $\cK$ plays a similar role for sequence algebras as the
ideal of the compact operators does for operator algebras. In
particular, there is a Fredholm theory for approximation sequences
which parallels the common Fredholm theory for operators and which
has remarkable consequences for the asymptotic behavior of
singular values (see Chapter 6 in \cite{HRS2} and Chapters 4 and 5
in \cite{Roc10} for an introduction and Section \ref{ss64} below
for a closer look at Fredholm properties of finite sections
sequences for operators in the Cuntz algebra).

The point now is that the sequences (\ref{e92.14}) and
(\ref{e92.15}) are compact in this sense (and do not belong to the
smaller ideal $\cG$ of the zero sequences). Since the Cuntz
algebras $\cO_N$ do not possess non-zero compact operators at all,
it seems to be not natural to consider discretizations which
produce non-zero compact sequences.

The second argument comes from operator theory. The following
lemma states that the choice (\ref{e92.16}) implies that, up to
sequences in the ideal $\cG_\eta$, the initial projections of the
partial isometries $(P_{\eta(n)} S_i P_{\eta(n)})$ with $i \in
\Omega^k$ only depend on the length $k$ of the multi-index, not on
the multi-index $i$ itself.
\begin{lemma} \label{l92.17}
Let $i \in \Omega^k$ and $n = N^j$ with $j \in \sZ^+$. Then
\[
P_n S_i^* P_n S_i P_n = \left\{
\begin{array}{ll}
0   & \mbox{if} \; j < k \; \mbox{and} \; N^j \le v_{i,k}, \\
P_1 & \mbox{if} \; j < k \; \mbox{and} \; N^j > v_{i,k}, \\
P_{N^{j-k}} & \mbox{if} \; j \ge k.
\end{array}
\right.
\]
\end{lemma}
{\bf Proof.} By (\ref{e92.11}), one has $P_n S_i^* P_n S_i P_n =
P_r$ with
\[
r = \{ (N^j - v_{i,k})/N^k \} = \{ N^{j-k} - v_{i,k}/N^k \}.
\]
If $j \ge k$, then $N^{j-k}$ is a positive integer, whereas
$v_{i,k}/N^k \in [0, \, 1)$. Thus, in this case, $r = N^{j-k}$.
Now let $j < k$. Then
\[
(N^j - v_{i,k})/N^k \le N^{j-k} < 1.
\]
Thus, if $N^j - v_{i,k} > 0$, then $r = 1$. Finally, let $j < k$
and $N^j - v_{i,k} \le 0$. Then we conclude from
\[
(N^j - v_{i,k})/N^k \ge N^{j-k} - 1 > -1
\]
that $r = 0$ (recall that we agreed upon $P_0 = 0$). \hfill \qed
\subsection{The restricted algebra $\cS_\eta (\cO_N)$} \label{ss34}
In what follows we will exclusively deal with the restricted
algebra $\cS_\eta(\cO_N)$ where $\eta(n) = N^n$, as suggested by
the arguments of the previous section.
\begin{prop} \label{p92.18}
The algebra $\cS_\eta (\cO_N)$ contains the ideal $\cG_\eta$.
\end{prop}
{\bf Proof.} In a first step we show that $\cS_\eta (\cO_N)$
contains all sequences of the form
\begin{equation} \label{e92.19}
(0, \, \ldots, \, 0, \, P_1, \, 0, \, 0, \, \ldots) \in \cF_\eta
\end{equation}
where $P_1$ stands at the $k$th position, $k$ an arbitrary
positive integer.

First let $k=1$. Define 1-tuples
\[
i_> := (0), \quad i_< := (1) \; \in \Omega^1.
\]
Then $v_{i_>,1} = 0$ and $v_{i_<,1} = 1$, and for $j=0$ one has
\[
N^j > v_{i_>,1} = 0 \quad \mbox{but} \quad N^j \le v_{i_<,1} = 1.
\]
By Lemma \ref{l92.17},
\[
(P_{N^n} S_{i_>}^* P_{N^n} S_{i_>} P_{N^n})_{n \ge 0} = (P_1, \,
P_{N^0}, \, P_{N^1}, \, P_{N^2}, \, \ldots) \in \cS_\eta(\cO_N)
\]
and
\[
(P_{N^n} S_{i_<}^* P_{N^n} S_{i_<} P_{N^n})_{n \ge 0} = (0, \,
P_{N^0}, \, P_{N^1}, \, P_{N^2}, \, \ldots) \in \cS_\eta(\cO_N).
\]
Hence, the sequence $(P_1, \, 0, \, 0, \, \ldots)$ belongs to
$\cS_\eta(\cO_N)$. Let now $k \ge 2$ and set
\[
i_> := (\underbrace{0, \, 0, \, \ldots, \, 0}_{k-2}, \, 1, \, 0), \quad
i_< := (\underbrace{0, \, 0, \, \ldots, \, 0}_{k-1}, \, 1) \in \Omega^k.
\]
Then $v_{i_>,k} = N^{k-2}$ and $v_{i_<,k} = N^{k-1}$, whence via
Lemma \ref{l92.17},
\[
(P_{N^n} S_{i_>}^* P_{N^n} S_{i_>} P_{N^n})_{n \ge 0} =
(\underbrace{0, \, 0, \, \ldots, \, 0}_{k-1}, \, P_1, \, P_{N^0},
\, P_{N^1}, \, P_{N^2}, \, \ldots) \in \cS_\eta(\cO_N)
\]
and
\[
(P_{N^n} S_{i_<}^* P_{N^n} S_{i_<} P_{N^n})_{n \ge 0} =
(\underbrace{0, \, 0, \, \ldots, \, 0}_{k}, \, P_{N^0}, \,
P_{N^1}, \, P_{N^2}, \, \ldots) \in \cS_\eta(\cO_N).
\]
Hence, the sequence $(0, \, \ldots, \, 0, \, P_1, \, 0, \,
\ldots)$ with $P_1$ standing at the $k$th position belongs to
$\cS_\eta(\cO_N)$, too.

In the next step we show that $\cS_\eta (\cO_N)$ contains all
sequences of the form
\begin{equation} \label{e92.20}
(0, \, \ldots, \, 0, \, A, \, 0, \, 0, \, \ldots) \in \cF_\eta
\end{equation}
where $A$ is an arbitrary $N^k \times N^k$ matrix standing at the
$k$th position. Since all sequences of the form (\ref{e92.19})
belong to $\cS_\eta (\cO_N)$, it is sufficient to check that the
set of all $n \times n$ matrices $P_1$, $P_n S_i P_n$ and $P_n
S_i^* P_n$ with $i \in \Omega^k$ and $k \ge 1$ generates the
algebra $\sC^{n \times n}$ (of course, we will need this fact
later on only for $n$ being a power of $N$). Let $0 \le r, \, s <
n$. Write $r$ and $s$ in the $N$-adic system as
\[
r =  i_1 + i_2N + \ldots + i_k N^{k-1}, \qquad
s =  j_1 + j_2N + \ldots + j_k N^{k-1}
\]
with $i = (i_1, \, \ldots, \, i_k), \; j  = (j_1, \, \ldots, \, j_k)
\in \Omega^k$. Then $r = v_{i,k}$ and $s = v_{j,k}$, and
\begin{eqnarray*}
P_n S_j P_n \cdot P_1 \cdot P_n S_i^* P_n (x_m)_{m = 0}^{n-1}
& = & (\underbrace{0, \, \ldots, \, 0}_{v_{j,k}}, \, x_{v_{i,k}}, \,
0, \, \ldots, \, 0) \\
& = & (\underbrace{0, \, \ldots, \, 0}_s, \, x_r, \,
0, \, \ldots, \, 0).
\end{eqnarray*}
Thus, $P_n S_j P_n \cdot P_1 \cdot P_n S_i^* P_n$ is a matrix
which takes the $r$th entry of a vector and writes it on the $s$th
place. Since every matrix is a linear combination of matrices of
this kind, the assertion of the second step follows.

Finally, every sequence in $\cG_\eta$ can be approximated as
closely as desired by finite sums of sequences of the form
(\ref{e92.20}). Since $\cS_\eta (\cO_N)$ is closed, this implies
$\cG_\eta \subseteq \cS_\eta (\cO_N)$. \hfill \qed
\section{The finite sections algebra $\cS_N$} \label{s4}
\subsection{A distinguished ideal} \label{ss41}
By Proposition \ref{p92.18}, one can form the quotient algebra
$\cS_\eta (\cO_N)/\cG_\eta$. We denote it by $\cS_N$. Recall that $\cS_N$ is generated by the partial isometries
\[
s_i := (P_{N^n} S_i P_{N^n})_{n \ge 0} + \cG_\eta, \quad i \in
\Omega
\]
and contains the identity element $e$ of $\cF_\eta/\cG_\eta$.
For each multi-index $i \in \Omega^k$ we further set
\[
s_i := s_{i_1} s_{i_2} \ldots s_{i_k}.
\]
These elements are partial isometries by Lemma \ref{l92.9},
and the initial projection of $s_i$ does only depend on the
length of $i$ by Lemma \ref{l92.17}. We denote the length of
the multi-index $i$ by $|i|$ and write $p_k$ for the joint
initial projection of all partial isometries $s_i$ with
length $k$. Further we write $\Omega_\infty$ for the set
of all multi-indices (of arbitrary length). The following
axioms collect the basic properties of these elements.
\begin{enumerate}
\item[(A1)]
for every $i \in \Omega_\infty$, the coset $s_i$ is a partial
isometry, the initial projection of which depends on $|i|$ only:
$s_i^* s_i = p_{|i|}$.
\item[(A2)]
$s_0 s_0^* + s_1 s_1^* + \ldots + s_{N-1} s_{N-1}^* = e$.
\end{enumerate}
All results in Section \ref{s4} will follow only from these
two axioms. For later reference we list some further relations between partial isometries $s_i$ and the projections $p_k$.
\begin{lemma} \label{l92.21}
Let $k, \, l$ positive integers and $i \in \Omega^k$. Then \\[1mm]
$(a)$ $s_i = s_i p_k$ and $s_i^* = p_k s_i^*$. \\
$(b)$ $s_i^* p_l s_i = p_{k+l}$. \\
$(c)$ $p_k p_l = p_k$ if $k \ge l$. \\
$(d)$ $p_l s_i = s_i p_{k+l}$. \\
$(e)$ The generalized Cuntz condition $\sum_{i \in \Omega^k}
s_i s_i^* = e$ holds for every $k \ge 1$.
\end{lemma}
{\bf Proof.} Assertion $(a)$ is the fact that $p_k$ is the initial
projection of $s_i$. For $(b)$ write $p_l$ as $s_j^* s_j$ with $j
\in \Omega^l$. Then $s_i^* p_l s_i = s_i^* s_j^* s_j s_i =
p_{k+l}$ by the definition of $p_{k+l}$. For $(c)$, write $p_k$
and $p_l$ as $(s_0^*)^k s_0^k$ and $(s_0^*)^l s_0^l$. Since
$s_0^l$ is a partial isometry,
\[
p_k p_l = (s_0^*)^k s_0^{k - l} s_0^l (s_0^*)^l s_0^l = (s_0^*)^k
s_0^{k - l} s_0^l = p_k,
\]
which gives $(c)$. Assertions $(b)$ and $(c)$ imply that
\begin{eqnarray*}
(p_l s_i - s_i p_{k+l})^* (p_l s_i - s_i p_{k+l}) & = & (s_i^* p_l
- p_{k+l} s_i^*) (p_l s_i - s_i p_{k+l}) \\
& = & s_i^* p_l s_i - s_i^* p_l s_i p_{k+l} - p_{k+l} s_i^* p_l
s_i + p_{k+l} s_i^* s_i p_{k+l} \\
& = & p_{k+l} - p_{k+l} - p_{k+l} + p_{k+l} p_k p_{k+l} = 0,
\end{eqnarray*}
whence $(d)$ via the $C^*$-axiom. Finally, assertion $(e)$ follows
easily by induction. For $k=1$, $(e)$ reduces to $(A2)$. If assertion
$(e)$ holds for some $k \ge 1$, then it holds for $k+1$ since
\[
\sum_{i \in \Omega^{k+1}} s_i s_i^*
= \sum_{(i_1, \tilde{i}) \in \Omega^{k+1}}
s_{(i_1, \tilde{i})} s_{(i_1, \tilde{i})}^*
= \sum_{i_1 \in \Omega^1} \sum_{\tilde{i} \in \Omega^k}
s_{i_1} s_{\tilde{i}} s_{\tilde{i}}^* s_{i_1}^*
= \sum_{i_1 \in \Omega^1} s_{i_1} s_{i_1}^* = e
\]
by assumption and axiom (A2). \hfill \qed \\[3mm]
For every positive integer $k$, let $\cJ^{(k)}$ denote the
smallest closed ideal of $\cS_N$ which contains the projection
$e - p_k$. By axiom (A1), every partial isometry $s_i$ of length
$k$ is an isometry modulo $\cJ^{(k)}$. By Lemma \ref{l92.17},
\[
e - p_1 = (0, \, P_N - P_1, \, P_{N^2} - P_N, \, P_{N^3} -
P_{N^2}, \, \ldots) + \cG_\eta.
\]
\begin{prop} \label{p92.22}
$\cJ^{(k)} = \cJ^{(1)}$ for every $k$.
\end{prop}
{\bf Proof.} By Lemma \ref{l92.21} $(c)$,
\[
(e-p_k)(e-p_1) = e - p_k - p_1 + p_k p_1 = e - p_1.
\]
Hence, $e-p_1 \in \cJ^{(k)}$, whence $\cJ^{(1)} \subseteq
\cJ^{(k)}$. For the reverse inclusion recall from Lemma
\ref{l92.21} $(b)$ that
\[
(s_0^*)^l (e-p_1) s_0^l = p_l - p_{l+1}
\]
for every $l \in \sZ^+$. Adding these identities for $l$ between 0
and $k-1$ gives $e - p_k$ on the right-hand side, whereas the
element of the left-hand side belongs to $\cJ^{(1)}$. Thus, $e-p_k
\in \cJ^{(1)}$, whence $\cJ^{(k)} \subseteq \cJ^{(1)}$. \hfill
\qed \\[3mm]
In what follows we write $\cJ_N$ for the ideal $\cJ^{(1)}$ of
$\cS_N$. Note that {\em every} partial isometry $s_i$ with
$i \in \Omega_\infty$ is an isometry modulo $\cJ_N$.
\begin{remark}
The smallest closed ideal of $\cS_N$ which contains the projection
$p_1$ coincides with all of $\cS_N$. Indeed, from $s_i = s_i s_i^*
s_i = s_i p_1$ we conclude that every element $s_i$ with $i \in
\Omega$ belongs to this ideal. Since the $s_i$ generate the
algebra $\cS_N$, the assertion follows.
\end{remark}
We conclude this section with a further property of the $p_n$
which will be needed in Section \ref{s4}.
\begin{lemma} \label{l92.23}
For each $j \in \cJ_N$, one has $\lim_{n \to \infty} \|p_n j\| =
0$.
\end{lemma}
{\bf Proof.} If $j$ is of the form
\begin{equation} \label{e92.24}
s_i s_k^* (e - p_1) s_l s_m^* \quad \mbox{with} \quad |i| \le |k|
\end{equation}
then the assertion holds since
\[
p_n s_i s_k^* (e - p_1) = s_i s_k^* p_{n + |i| - |k|} (e - p_1) =
0
\]
for $n > |k| - |i|$. Hence, the assertion also holds if $j$ is a
linear combination of elements of the form (\ref{e92.24}). Since
these linear combinations form a dense subset of $\cJ_N$ and
$\|p_n\| = 1$ for each $n$, the assertion holds for every $j \in
\cJ_N$. \hfill \qed
\subsection{Lifting $\cS_N/\cJ_N$ to $\cS_N$} \label{ss42}
Our further analysis of the algebra $\cS_N$ is based on the
following elementary fact, which can be considered as the simplest
instance of a lifting theorem. It settles a condition which
guarantees than every element which is invertible modulo an ideal
can be lifted to an invertible element. The simple proof is in
\cite{HRS2}, Theorem 5.33.
\begin{prop} \label{p93.1}
Let $\cA$ be a unital $C^*$-algebra and $\cI$ a closed ideal of
$\cA$. Further suppose there is a unital $^*$-homomorphism $\pi$
from $\cA$ into a unital $C^*$-algebra $\cB$ such that the
restriction of $\pi$ onto $\cI$ is injective. Then the following
assertions are equivalent for every element $a \in \cA$: \\[1mm]
$(a)$ $a$ is invertible in $\cA$. \\
$(b)$ The coset $a + \cI$ is invertible in the quotient algebra
$\cA/\cI$, and $\pi(a)$ is invertible in $\cB$.
\end{prop}
We shall apply this result with $\cA := \cS_N$ and $\cI := \cJ_N$.
By Proposition \ref{p93.1}, the problem to derive a criterion for
the invertibility of elements of $\cS_N$ (and thus, for the
stability of sequences in $\cS_\eta(\cO_N)$) splits into two
separate tasks:
\begin{itemize}
\item
to describe the quotient algebra $\cS_N/\cJ_N$, and
\item
to construct an injective $^*$-homomorphism on $\cJ_N$.
\end{itemize}
The solution of the first task is evident: The quotient algebra
$\cS_N/\cJ_N$ is generated by the cosets $s_i + \cJ_N$ with
$i \in \Omega$. These cosets are isometries and they satisfy
the Cuntz axiom
\[
(s_0 + \cJ_N) (s_0 + \cJ_N)^* + \ldots + (s_{N-1} + \cJ_N)
(s_{N-1} + \cJ_N)^* = e + \cJ_N.
\]
By the universal property of Cuntz algebras, $\cS_N/\cJ_N$
is $^*$-isomorphic to the (abstract) Cuntz algebra $\cO_N$.
It is also not difficult to construct an isomorphism from
$\cS_N/\cJ_N$ onto the (concrete) Cuntz algebra $\cO_N$
explicitly. For, let $W_\eta : \cF_\eta \to L(l^2(\sZ^+))$
denote the mapping which associates with each sequence in
$\cF_\eta$ its strong limit. Clearly, $W_\eta$ is a
$^*$-homomorphism. Since the ideal $\cG_\eta$ lies in the
kernel of $W_\eta$, there is a correctly defined quotient
homomorphism
\begin{equation} \label{e93.2}
W_\eta^\cG : \cF_\eta/\cG_\eta \to L(l^2(\sZ^+)), \quad \bA +
\cG_\eta \mapsto W_\eta (\bA).
\end{equation}
Applying this homomorphism to both sides of the equality $s_0^*
s_0 = p_1$ we get $S_0^* S_0 = W_\eta^\cG (p_1)$, whence
$W_\eta^\cG (p_1) = I$. Hence, the ideal $\cJ_N$ lies in the
kernel of $W_\eta^\cG$, which implies that the quotient homomorphism
\begin{equation} \label{e93.3}
(\cS_\eta(\cO_N)/\cG_\eta)/\cJ_N \to L(l^2(\sZ^+)), \quad
(\bA + \cG_\eta) + \cJ_N \mapsto W_\eta^\cG (\bA + \cG)
\end{equation}
is correctly defined, too; we denote it by $W^\cJ$.
\begin{theo} \label{t93.4}
$W^\cJ$ is a $^*$-isomorphism from $\cS_N/\cJ_N$ onto the
(concrete) Cuntz algebra $\cO_N$.
\end{theo}
{\bf Proof.} The $^*$-homomorphism $W^\cJ$ maps the generating
cosets $s_i + \cJ_N$, $i \in \Omega$, to the generating operators
$S_i$ of the Cuntz algebra $\cO_N$, respectively. Since both sets
of generators consist of partial isometries which satisfy the
(same) Cuntz axiom, the assertion follows form the universal
property of Cuntz algebras again. \hfill \qed \\[3mm]
The following is an immediate consequence of this theorem and of
the fact that $\cO_N$ is a simple algebra.
\begin{coro} \label{c93.5}
The kernel of the restriction of the homomorphism $W_\eta^\cG$
defined by $(\ref{e93.2})$ to the algebra $\cS_N$ coincides with
$\cJ_N$.
\end{coro}
The following fact sheds a first light on our second task. Some
consequences of this fact are already discussed in the
introduction.
\begin{theo} \label{t93.6}
Every proper closed ideal of $\cS_N$ lies in $\cJ_N$.
\end{theo}
{\bf Proof.} Let $\widetilde{\cJ}$ be a proper closed ideal of
$\cS_N$. Then $\cJ_N + \widetilde{\cJ}$ is a closed ideal of
$\cS_N$ with $\cJ_N \subseteq \cJ_N + \widetilde{\cJ} \subseteq
\cS_N$. Since the quotient $\cS_N/\cJ_N$ is $^*$-isomorphic
to $\cO_N$ and, hence, a simple algebra, one has either
\begin{itemize}
\item
Case A: $\cJ_N + \widetilde{\cJ} = \cS_N$, or
\item
Case B: $\cJ_N + \widetilde{\cJ} = \cJ_N$, i.e. $\widetilde{\cJ}
\subseteq \cJ_N$.
\end{itemize}
We wish to exclude case A. Suppose we are in the situation of
case A. Consider the ideals $\cI_1 := \cJ_N/(\cJ_N \cap
\widetilde{\cJ})$ and $\cI_2 := \widetilde{\cJ}/(\cJ_N \cap
\widetilde{\cJ})$ of $\cB := \cS_N/(\cJ_N \cap \widetilde{\cJ})$.
These ideals have a trivial intersection, their sum is $\cB$,
and the algebra
\[
\cB/\cI_1 = \left( \cS_N/(\cJ_N \cap \widetilde{\cJ}) \right) /
\left( \cJ_N/(\cJ_N \cap \widetilde{\cJ}) \right) \cong
\cS_N/\cJ_N
\]
is still simple. Let $W$ stand for the canonical homomorphism
from $\cB$ onto $\cB/\cI_2$ and write $\hat{a}$ for the coset
of $a \in \cS_N$ modulo $\cJ_N \cap \widetilde{\cJ}$. Since
$W(\cB) = W(\cI_1)$, there is an element $\hat{\pi} \in \cI_1$
such that $W(\hat{\pi}) = W(\hat{e})$. From
\[
W (\hat{\pi}^2 - \hat{\pi}) = W (\hat{e}^2 - \hat{e}) = 0 \quad
\mbox{and} \quad W (\hat{\pi}^* - \hat{\pi}) = W (\hat{e}^*
- \hat{e}) = 0
\]
we conclude that both $\hat{\pi}^2 - \hat{\pi}$ and $\hat{\pi}^*
- \hat{\pi}$ belong to $\cI_1 \cap \cI_2$. Since this intersection
is trivial, the element $\hat{\pi}$ is a (self-adjoint) projection.
Moreover, since
\[
W (\hat{a} \hat{\pi} - \hat{\pi} \hat{a}) = W (\hat{a}) W
(\hat{e}) - W (\hat{e}) W (\hat{a}) = 0
\]
for every element $\hat{a} \in \cB$ we conclude as above that
$\hat{\pi}$ lies in the commutant of $\cB$. A similar reasoning
shows finally that $\hat{\pi}$ is the identity element for
$\cI_1$. Similarly, $\hat{e} - \hat{\pi}$ belongs to $\cI_2$ and
is the identity element for $\cI_2$.

Let $\pi \in \cJ_N$ be a representative of the coset $\hat{\pi}$.
From Lemma \ref{l92.23} we infer
\[
\|(e - p_n) \pi - \pi\| \to 0 \quad \mbox{as} \quad n \to \infty,
\]
whence
\[
\|(\hat{e} - \widehat{p_n}) \hat{\pi} - \hat{\pi}\| \to 0 \quad
\mbox{as} \quad n \to \infty.
\]
Since $\hat{\pi}$ is the identity element of $\cI_1$ and $\hat{e}
- \widehat{p_n} \in \cI_1$, this implies
\[
\|\hat{e} - \widehat{p_n} - \hat{\pi}\| \to 0 \quad
\mbox{as} \quad n \to \infty.
\]
Since $\hat{e} - \widehat{p_n}$ and $\hat{\pi}$ are commuting
projections, this finally shows that $\hat{e} - \widehat{p_n}
= \hat{\pi}$ for all sufficiently large $n$, say $n \ge k$.
Consequently, for $n \ge k$, one has $\widehat{p_n} = \hat{e}
- \hat{\pi} \in \cI_2$, whence $p_n \in \widetilde{\cJ}$. Since
$s_i = s_i p_k$ for all $i \in \Omega^k$ by Lemma \ref{l92.21},
this implies $s_i \in \widetilde{\cJ}$ and, thus, the smallest
closed ideal of $\cS_N$ which contains all partial isometries
$s_i$ with $|i| = k$ lies in $\widetilde{\cJ}$. By Lemma
\ref{l92.21} $(e)$, this finally implies $e \in \widetilde{\cJ}$.
Thus, $\widetilde{\cJ}$ is not a proper ideal of $\widetilde{\cJ}$.
This contradiction excludes case A. \hfill \qed
\begin{coro} \label{c93.8}
Every $^*$-homomorphism on $\cS_N$ which is injective on $\cJ_N$
is injective on all of $\cS_N$.
\end{coro}
Indeed, if $W$ is a $^*$-homomorphism on $\cS_N$ which is injective
on $\cJ_N$, then its kernel is a proper ideal of $\cS_N$. By Theorem
\ref{t93.6}, $\ker W \subset \cJ_N$. But $\cJ_N \cap \ker W = \{0\}$
by assumption. Hence, the kernel of $W$ is trivial, and $W$ is
injective on $\cS_N$. \hfill \qed \\[3mm]
In Section \ref{s6} we are going to construct an injective
homomorphism on $\cJ_N$. We prepare this construction by a closer
look at the Cuntz algebra and a related Toeplitz algebra in
Section \ref{s5}
\section{Expectations of $\cO_N$ and stratified Toeplitz
operators} \label{s5}
\subsection{Block Toeplitz operators} \label{ss51}
By $l^2(\sZ^+, \, l^2(\sZ^+))$ we denote the Hilbert space of all
sequences $x = (x_n)_{n \ge 0}$ with values in $l^2(\sZ^+)$ such
that
\[
\|x\|^2 := \sum_{n \ge 0} \|x_n\|^2 < \infty.
\]
It will be convenient to identify the algebra of all bounded
linear operators on $l^2(\sZ^+, \, l^2(\sZ^+))$ with the minimal
tensor product $L(l^2(\sZ^+)) \otimes L(l^2(\sZ^+))$. Let $A$ be
a bounded linear operator on $l^2(\sZ^+, \, l^2(\sZ^+))$, and let
$A_{ij} \in L(l^2(\sZ^+))$ be the operator which maps the $i$th
component of $x$ to the $j$th component of $Ax$. Then $A$ can
be identified with the infinite matrix $(A_{ij})_{i, j \ge 0}$
in an evident way (but, of course, not every matrix with entries
in $L(l^2(\sZ^+))$ defines a bounded operator on $l^2(\sZ^+, \,
l^2(\sZ^+))$).

The closure in $L(l^2(\sZ^+, \, l^2(\sZ^+)))$ of the set of all
matrices with only finitely non-vanishing entries forms a closed
ideal of this algebra which we identify with $K(l^2(\sZ^+))
\otimes L(l^2(\sZ^+))$ in a natural way.

To each function $a \in C(\sT, L(l^2(\sZ^+)))$ we associate its
$k$th Fourier coefficient
\[
a_k := \int_\sT a(\lambda) \lambda^{-k} \, d\lambda \in L(l^2(\sZ^+)), \quad k \in \sZ,
\]
and consider the (infinite) Toeplitz matrix
\[
T(a) := (a_{i-j})_{i,j \ge 0}.
\]
Every Toeplitz matrix with a continuous generating function
defines a bounded operator on $l^2(\sZ^+, \, l^2(\sZ^+))$
which is a called a Toeplitz operator and denoted by $T(a)$
again, and
\begin{equation} \label{e94.1}
\|T(a)\| = \|a\|_\infty.
\end{equation}
The Toeplitz operator $T(a)$ is invertible modulo $K(l^2(\sZ^+))
\otimes L(l^2(\sZ^+))$ if and only if the function $a$ is
invertible in $C(\sT, \, L(l^2(\sZ^+)))$. For details see
\cite{Pag1}.
\subsection{An algebra of stratified Toeplitz operators}
\label{ss52}
For $i \in \Omega$, consider the infinite matrix
\begin{equation} \label{e94.2}
\Sigma_i := \pmatrix{0 & S_i &     &        & \cr
                       & 0   & S_i &        & \cr
                       &     & 0   & S_i    & \cr
                       &     &     & \ddots & \ddots}
\end{equation}
with all entries left empty being zero, and for each multi-index
$i \in \Omega^k$, let
\[
\Sigma_i := \Sigma_{i_1} \ldots \Sigma_{i_k}.
\]
Clearly, every $\Sigma_i$ is a Toeplitz matrix with continuous 
generating function. It defines a bounded operator on
$l^2(\sZ^+, \, l^2(\sZ^+))$, which we denote by $\Sigma_i$ again.
We let $\cT_N$ refer to the smallest closed subalgebra of
$L(l^2(\sZ^+, \, l^2(\sZ^+)))$ which contains all operators
$\Sigma_i$ and $\Sigma_i^*$ with $i \in \Omega$.

One peculiarity of operators in $\cT_N$ which will become
more clear and important later on should be mentioned
already here: Their matrix representation is stratified in
the sense that all entries on the $k$th diagonal above the
main diagonal are necessarily of the form $\Sigma_i \cO_N^{par}$
with an multi-index $i$ of length $k$. Thus, the only operator
which can stand on each diagonal is the zero operator. (A
curious consequence of this fact is that the only block Hankel operators which are contained in the Toeplitz algebra $\cT_N$
are of the form of a diagonal matrix $\diag (A, \, 0, \, 0, \, \ldots)$ with $A \in L(l^2(\sZ^+))$.)

One easily checks that the operators $\Sigma_i$ are partial
isometries which satisfy the Cuntz axiom
\[
\sum_{i \in \Omega} \Sigma_i \Sigma_i^* = \sum_{i \in \Omega}
\diag(S_i S_i^*, \, S_i S_i^*, \ldots) = \diag (I, \, I, \,
\ldots).
\]
Moreover, for each multi-index $i$, the initial projection
$\Sigma_i^* \Sigma_i$ is equal to $I - \Pi_{|i|}$ where
\begin{equation} \label{e94.3}
\Pi_k := \diag (\underbrace{I, \, \ldots, I}_k, \, 0, \, 0, \,
\ldots) \in \cT_N
\end{equation}
for $k \ge 1$. Thus, the algebra $\cT_N$ contains the identity
operator $I$ and all projections $\Pi_k$, and the partial
isometries $\Sigma_i$ satisfy the axioms (A1) and (A2)
in Section \ref{ss41} in place of the $s_i$. Thus, all results
of Section \ref{s4} will remain valid for the algebra $\cT_N$
in place of $\cS_N$ and for its ideal $\cC_N$, which is the
smallest closed ideal of $\cT_N$ which contains the projection
$\Pi_1$, in place of $\cJ_N$. In particular,
\begin{equation} \label{e94.4}
\Pi_n \in \cC_N \quad \mbox{for each} \; n \ge 1
\end{equation}
and
\begin{equation} \label{e94.5}
\cT_N/\cC_N \cong \cO_N.
\end{equation}
\subsection{Identification of $\cC_N$} \label{ss53}
Our next goal is a description of the ideal $\cC_N$ of $\cT_N$.
Let $\cO_N^{par}$ refer to the smallest closed subalgebra of
$\cO_N$ which contains all products $S_i S_j^*$ with multi-indices $i$ and $j$ of the same length $|i| = |j|$. Here we allow
multi-indices of length 0 and set $S_{\emptyset} := I$. Thus,
$\cO_N^{par}$ is a unital algebra. One easily checks that $\cO_N$ is the closed span of the set of all products $S_i S_j^*$, whereas $\cO_N^{par}$ is the closed span of all products $S_i S_j^*$ with $|i| = |j|$. The algebra $\cO_N^{par}$ is known to be isomorphic to the UHF-algebra of type $N^\infty$. In particular, $\cO_2^{par}$ is the CAR-algebra. See \cite{Cun1,Dav1} for details.

We shall further make use of the following elementary observation.
Let $\cA$ be a $C^*$-subalgebra of a unital $C^*$-algebra $\cB$
and let $\Sigma$ be an isometry in $\cB$. Then the mapping
\begin{equation} \label{e94.6}
\cA \to \cB, \quad A \mapsto \Sigma A \Sigma^*
\end{equation}
is an injective $^*$-homomorphism. Indeed, $\Sigma A_1 \Sigma^*
\Sigma A_2 \Sigma^* = \Sigma A_1 A_2 \Sigma^*$ for $A_1, \, A_2
\in \cA$, and if $\Sigma A \Sigma^* = 0$, then
\[
0 = \Sigma^* \Sigma A \Sigma^* \Sigma = A.
\]
We apply this observation to the algebras $\cA := \cC_N$ and $\cB
:= L(l^2(\sZ^+, \, l^2(\sZ^+)))$ and to the isometry
\[
\Sigma := \diag (I, \, S_0, \, S_0^2, \, S_0^3, \, \ldots).
\]
Thus, $\cC_N$ is $^*$-isomorphic to $\Sigma \cC_N \Sigma^*$. The
latter algebra can be described as follows where we let $\Pi :=
\Sigma \Sigma^*$ be the image of the identity operator under the
mapping (\ref{e94.6}).
\begin{theo} \label{t94.7}
$\Sigma \cC_N \Sigma^* = \Pi \left( K(l^2(\sZ^+)) \otimes
\cO_N^{par} \right) \Pi$.
\end{theo}
{\bf Proof.} Since the algebra $\cT_N$ is generated by the
partial isometries $\Sigma_i$, the algebra $\Sigma \cT_N \Sigma^*$
is generated by the operators
\[
\Sigma \Sigma_i \Sigma^* = \pmatrix{0 & S_i S_0^* &  &   & \cr
                       & 0   & S_0 S_i (S_0^*)^2 &        & \cr
                       &     & 0   & S_0^2 S_i (S_0^*)^3  & \cr
                       &     &     & \ddots & \ddots}.
\]
All entries of this matrix belong to $\cO_N^{par}$. Hence,
\[
\Sigma \cT_N \Sigma^* \subseteq \Pi \left( L(l^2(\sZ^+)) \otimes
\cO_N^{par} \right) \Pi.
\]
Further, the mapping (\ref{e94.6}) sends the generator $\Pi_1$
of the ideal $\cC_N$ to itself. Since $\Pi_1 \Pi = \Pi \Pi_1$,
one has
\[
\Sigma \Pi_1 \Sigma^* = \Pi_1 \in \Pi \left( K(l^2(\sZ^+)) \otimes
\cO_N^{par} \right) \Pi,
\]
whence the inclusion
\[
\Sigma \cC_N \Sigma^* \subseteq \Pi \left( K(l^2(\sZ^+)) \otimes
\cO_N^{par} \right) \Pi.
\]
The reverse inclusion will follow once we have shown that for
each $n \times n$-matrix $A := (A_{ij})$ with entries in
$\cO_N^{par}$, which we identify with an operator on the range of
$\Pi_n$ in the obvious way, the operator $\Pi A \Pi$ belongs to
$\Sigma \cC_N \Sigma^*$. Due to linearity we can further assume
that only one of the entries of $A$, say $A_{ij}$, is different
from 0. Finally, since $\cO_N^{par}$ is spanned by products
$S_l S_r^*$ with multi-indices of the same length, we can
assume that only the $ij$th entry of $A$ is different from zero
and that this entry is $S_l S_r^*$ with $|r| = |l|$. Then
$\Pi A \Pi$ is a matrix the only non-vanishing entry of which
stands at the $ij$th position, and this entry is
\begin{equation} \label{e94.8}
S_0^i (S_0^*)^i \, S_l S_r^* \, S_0^j (S_0^*)^j.
\end{equation}
Let
\[
B := (\Pi_{i+1} - \Pi_i) (\Sigma_0^*)^i \Sigma_l \Sigma_r^*
\Sigma_0^j (\Pi_{j+1} - \Pi_j).
\]
This operator is in $\cC_N$, all entries in the matrix
representation of $\Sigma B \Sigma^*$ with exception of
the $ij$th entry vanish, and the $ij$th entry coincides with
(\ref{e94.8}). Thus, $\Sigma B \Sigma^* = \Pi A \Pi$, which
finishes the proof.
\hfill \qed \\[3mm]
Next we will have a closer look at the ideal structure of
$\cT_N$.
\begin{theo} \label{t94.9}
The ideal $\cC_N$ of $\cT_N$ is simple.
\end{theo}
{\bf Proof.} Let $\cR$ be a closed ideal of $\cC_N$. Then
$\Pi_1 \cR \Pi_1$ is a closed ideal of $\Pi_1 \cC_N \Pi_1$.
From Theorem \ref{t94.7} we infer that $\Pi_1 \cC_N \Pi_1$
is $^*$-isomorphic to the algebra $\cO_N^{par}$ which on
its hand is known to be $^*$-isomorphic to the UHF-algebra of
type $N^\infty$. Thus, $\Pi_1 \cC_N \Pi_1$ is isomorphic to
the inductive limit
\[
\sC \to \sC^{N \times N} \to \sC^{N^2 \times N^2} \to
\sC^{N^3 \times N^3} \to \ldots
\]
with connecting maps
\[
a \mapsto \diag (\underbrace{a, \, a, \, \ldots, a}_N).
\]
Being an inductive limit of simple algebras, the algebra
$\Pi_1 \cC_N \Pi_1$ is simple. Hence,
\[
\mbox{either} \quad \Pi_1 \cR \Pi_1 = \Pi_1 \cC_N \Pi_1
\quad \mbox{or} \quad \Pi_1 \cR \Pi_1 = \{0\}.
\]
In the first case, $\Pi_1 \in \Pi_1 \cR \Pi_1 \subseteq
\cR$. Since $\Pi_1$ generates $\cC_N$ as an ideal, we
conclude $\cR = \cC_N$.

Assume now that $\Pi_1 \cR \Pi_1 = \{0\}$. Let $R =
(R_{ij})_{ij \ge 0} \in \cR$. Then, for arbitrary subscripts
$i_0, \, j_0 \ge 0$, the matrix
\[
(\Pi_{i_0 + 1} - \Pi_{i_0}) R (\Pi_{j_0 + 1} - \Pi_{j_0})
\]
has the entry $R_{i_0j_0}$ at the $i_0j_0$th position
whereas all other entries are zero. Let $k$ and $l$ be
multi-indices with $|k| = i_0$ and $|l| = j_0$. Then the
matrix
\begin{equation} \label{e94.10}
\Sigma_k (\Pi_{i_0 + 1} - \Pi_{i_0}) R (\Pi_{j_0 + 1}
- \Pi_{j_0}) \Sigma_l^*
\end{equation}
has the entry $S_k R_{i_0j_0} S_l^*$ at the $00$th
position whereas all other entries are zero. Thus, the
matrix (\ref{e94.10}) belongs to $\Pi_1 \cR \Pi_1$, whence
$S_k R_{i_0j_0} S_l^* = 0$ by assumption. Since the $S_i$
are isometries, this implies $R_{i_0j_0} = 0$, and since
$i_0$ and $j_0$ were arbitrarily chosen, $R$ is the zero
matrix. Thus, $\cR$ is the zero ideal in this case. \hfill
\qed \\[3mm]
As already mentioned, the $\Sigma_i$ satisfy the axioms (A1)
and (A2) and, thus, all results of Section \ref{s4} hold
for the algebra $\cT_N$ in place of $\cS_N$ as well. In
particular, every proper closed ideal of $\cT_N$ lies in
$\cC_N$ by Theorem \ref{t93.6}. Together with Theorem
\ref{t94.9} this implies
\begin{coro} \label{c94.11}
$\cC_N$ is the only non-trivial closed ideal of $\cT_N$.
\end{coro}
\subsection{Expectations on $\cO_N$ and Toeplitz operators}
\label{ss54}
There are at least two ways to associate with every element
of the Cuntz algebra $\cO_N$ a Toeplitz operator in $\cT_N$ .
For facts cited without proof see \cite{Cun1,Dav1}.

The first way is via a special expectation. Recall that the
operators $S_i$ and $S_j^*$ with $i, j \in \Omega$ generate a
dense subalgebra of $\cO_N$. Each operator $A$ in this algebra can
be uniquely written as a finite sum
\begin{equation} \label{e94.12}
A = \sum_{k < 0} (S_0^*)^{-k} A_k + A_0 + \sum_{k > 0} A_k S_0^k
\end{equation}
with "Fourier" coefficients $A_k \in \cO_N^{par}$. For $k \in \sZ$
and $A$ as in (\ref{e94.12}), define $\Phi_k(A) := A_k$. Then
$\|\Phi_k (A)\| \le \|A\|$, thus the $\Phi_k$ extend by continuity
to bounded mappings from $\cO_N$ onto $\cO_N^{par}$. These
mappings own the following properties:
\begin{itemize}
\item
$\Phi_0 : \cO_N \to \cO_N^{par}$ is an expectation, i.e. $\Phi_0^2
= \Phi_0$,
\item
$\Phi_{k+1} (A) = \Phi_k(A S_0^*)$ if $k \ge 0$, and
\item
$\Phi_{k-1} (A) = \Phi_k(S_0 A)$ if $k < 0$.
\end{itemize}
We associate with each operator $A \in \cO_N$ a formal matrix of
operators on $l^2(\sZ^+)$ by
\begin{equation} \label{e94.13}
\Psi(A) := \left( (S_0^*)^i \Phi_0(S_0^i A (S_0^*)^j) S_0^j
\right)_{i,j \ge 0}.
\end{equation}
We will see in a moment that the formal matrix $\Psi(A)$ defines
a bounded operator on $l^2(\sZ^+, \, l^2(\sZ^+))$ and that this
operator is a Toeplitz operator in $\cT_N$. The following example
shows that this is at least true for operators in a dense
subalgebra of $\cO_N$.
\begin{example} \label{ex94.14}
Let $A := S_l S_m^*$ with multi-indices $l$ and $m$. Then
$\Phi_0(S_0^i A (S_0^*)^j)$ is different from 0 only if $i + |l| =
|m| + j$. In this case,
\[
\Phi_0(S_0^i A (S_0^*)^j) = \Phi_0(S_0^i S_l S_m^* (S_0^*)^j) =
S_0^i S_l S_m^* (S_0^*)^j,
\]
whence
\[
(S_0^*)^i \Phi_0(S_0^i A (S_0^*)^j) S_0^j = S_l S_m^* = A.
\]
Thus, $\Psi(S_l S_m^*)$ is the matrix which has the entry $S_l
S_m^*$ on its $i-j = |m| - |l|$th diagonal whereas all other
entries are zero. In particular, $\Psi(S_l S_m^*)$ is a Toeplitz
operator in $\cT_N$.

Note that the mapping $\Psi$ is not multiplicative. Indeed,
$\Psi(S_k) = \Sigma_k$ for $k \in \Omega$ as we have just seen.
Consequently, $\Psi(S_k^*) \Psi(S_k) = \Sigma_k^* \Sigma_k =
I - \Pi_1$, which is different from $\Psi(S_k^* S_k) =
\Psi(I) = I$. \hfill \qed
\end{example}
A second way to associate with every operator in $\cO_N$ a
Toeplitz operator in $\cT_N$ is via continuous functions. Let
$\lambda \in \sT$. Then the mapping $\rho_\lambda: S_i \mapsto
\bar{\lambda} S_i$ extends to an automorphism of $\cO_N$. Here, as
usual, $\bar{\lambda}$ stands for the complex conjugate of
$\lambda$; note that the mapping $\rho_\lambda$ is defined in
\cite{Dav1} without the bar. For each operator $A \in \cO_N$,
consider the function
\begin{equation} \label{e94.15}
f_A : \sT \to \cO_N, \quad \lambda \mapsto \rho_\lambda (A).
\end{equation}
\begin{lemma} \label{l94.16}
The function $f_A$ is continuous for each $A \in \cO_N$, and
$\|f_A\|_\infty = \|A\|$.
\end{lemma}
{\bf Proof.} For each $\lambda \in \sT$, one has $\|f_A(\lambda)\|
= \|\rho_\lambda (A)\| \le \|A\|$, whence $\|f_A\|_\infty \le
\|A\|$. Since $f_A(1) = A$, equality holds in this estimate.
Since $A \mapsto f_A$ is a linear mapping, this implies
\begin{equation} \label{e94.17}
\|f_A - f_B\|_\infty = \|A - B\| \quad \mbox{for all} \; A, \, B
\in \cO_N.
\end{equation}
Choose operators $B_n$ in the dense subalgebra of $\cO_N$
generated by the isometries $S_i$ such that $\|A - B_n\| \to 0$ as
$n \to \infty$. The functions $f_{B_n}$ are evidently continuous.
Being a uniform limit of continuous functions by (\ref{e94.17}),
the function $f_A$ is continuous. \hfill \qed \\[3mm]
As in Section \ref{ss51}, we associate with the continuous
function $f_A$ the sequence of its Fourier coefficients and
consider the associated Toeplitz operator $T(f_A)$ on
$l^2(\sZ^+, l^2(\sZ^+))$. From (\ref{e94.1}) and Lemma
\ref{l94.16} we conclude that
\begin{equation} \label{e94.18}
\|T(f_A)\| = \|f_A\|_\infty = \|A\| \quad \mbox{for every} \;
A \in \cO_N.
\end{equation}
\begin{example} \label{ex94.19}
Let $A = S_l S_m^*$ with multi-indices $l$ and $m$. Then
\[
f_A(\lambda) = \lambda^{-|l|} S_l \lambda^{|m|} =
\lambda^{|m|-|l|} A.
\]
The $|m|-|l|$th Fourier coefficient of $f_A$ is $A$ whereas all
other Fourier coefficients of this function vanish. Thus, $T(f_A)$
is the matrix which has the entry $S_l S_m^*$ on its $i-j = |m|
- |l|$th diagonal whereas all other entries are zero. \hfill \qed
\end{example}
Since the products $S_l S_m^*$ span a dense subalgebra of
$\cO_N$ and since the mapping $A \mapsto T(f_A)$ is an isometry,
we conclude from this example that $T(f_A)$ is a Toeplitz
operator {\em in} $\cT_N$ for every $A \in \cO_N$.

We will see now that the two ways discussed above lead to the
same goal.
\begin{theo} \label{t94.20}
The mapping $\Psi$ is a linear contraction from $\cO_N$ into
$\cT_N$. It coincides with the mapping $A \mapsto T(f_A)$.
\end{theo}
{\bf Proof.} We learned from Examples \ref{ex94.14} and
\ref{ex94.19} that $T(f_A) = \Psi(A)$ for $A = S_l S_m^*$. Since
these products span a dense subalgebra of $\cO_N$ and $A \mapsto
T(f_A)$ and $\Psi$ are linear mappings, this implies that
\begin{equation} \label{e94.21}
T(f_A) = \Psi(A) \in \cT_N \quad \mbox{for all} \; A \; \mbox{in a
dense subalgebra of} \; \cO_N.
\end{equation}
Since the mapping $A \mapsto T(f_A)$ is an isometry, we get from
(\ref{e94.21}) that
\[
\|\Psi(A)\| = \|T(f_A)\| = \|A\|
\]
for all $A$ in a dense subalgebra of $\cO_N$. Thus, the mapping $\Psi$ can be continued to a linear contraction from
$\cO_N$ into $\cT_N$ (which, of course, coincides with the mapping
$A \mapsto T(f_A)$). Since each entry of the matrix (\ref{e94.13})
depends continuously on $A$, this contractive continuation
coincides with the formal matrix in (\ref{e94.13}). \hfill \qed
\\[3mm]
The classical Toeplitz algebra decomposes into the direct sum of
the linear space $\{T(f): f \in C(\sT)\}$ and the ideal of the
compact operators. A similar decomposition holds for the algebra
$\cT_N$.
\begin{theo} \label{t94.22}
$\cT_N = \{ T(f_A): A \in \cO_N\} \oplus \cC_N = \{ \Psi(A): A \in
\cO_N\} \oplus \cC_N$.
\end{theo}
{\bf Proof.} We reify Proposition \ref{p91.7} with the following
algebras and mappings:
\begin{itemize}
\item
$\cA$ is the smallest closed subalgebra of $L(l^2(\sZ, \,
l^2(\sZ^+)))$ which contains all operators of Laurent type
represented by the two-sided infinite matrix
\[
\pmatrix{\ddots & \ddots &     &     &        & \cr
                & 0      & S_i &     &        & \cr
                &        & 0   & S_i &        & \cr
                &        &     & 0   & S_i    & \cr
                &        &     &     & \ddots & \ddots}
\]
with the zeros standing on the main diagonal.
\item
$\cB$ is the algebra $\cT_N$.
\item
$D$ is the mapping $\cA \to \cB, \; A \mapsto PAP$ where $P$ is the orthogonal projection from $l^2(\sZ, \, l^2(\sZ^+))$ onto $l^2(\sZ^+, \, l^2(\sZ^+))$ (the latter is identified with a closed subspace of the former in the obvious way).
\item
$W$ is the mapping
\[
\cB \to \cA, \quad B \mapsto \mbox{s-lim}_{n \to + \infty} V_n^* B V_n
\]
where $V$ is the operator of forward shift on $l^2(\sZ, \,
l^2(\sZ^+))$ and $V_n := V^n$ for $n \ge 0$, and where
$\mbox{s-lim}$ refers to the limit in the strong operator
topology.
\end{itemize}
Then Proposition \ref{p91.7} implies that $\cT_N = \{ T(f_A): A
\in \cO_N\} \oplus \ker W$, and it remains to verify that
\begin{equation} \label{e94.23}
\cC_N \; (:= \mbox{closid}_{\cT_N} \, \{ \Pi_1 \}) = \ker W.
\end{equation}
Evidently, $\Pi_1 \in \ker W$, whence the inclusion $\cC_N
\subseteq \ker W$. To get the reverse implication, we show that
the quasicommutator ideal of $\cT_N$ lies in $\cC_N$. Since the
products $S_l S_m^*$ with multi-indices $l, \, m$ span a dense
subalgebra of $\cO_N$, this fact will follow once we have shown
that
\begin{equation} \label{e94.24}
\Psi(S_l S_m^*) \Psi(S_n S_r^*) - \Psi(S_l S_m^* S_n S_r^*) \in
\cC_N
\end{equation}
for each choice of multi-indices $l, \, m, \, n$ and $r$. Let
$\Lambda$ denote the operator of forward shift
\begin{equation} \label{e94.25}
\Lambda := \pmatrix{0 &   &        & \cr
                    I & 0 &        & \cr
                      & I & 0      & \cr
                      &   & \ddots & \ddots}
\in L(l^2(\sZ^+, l^2(\sZ^+)))
\end{equation}
and set $\Lambda_n := \Lambda^n$ for $n > 0$, $\Lambda_0 := I$ and
$\Lambda_n := (\Lambda^*)^{-n}$ for $n < 0$. In Example \ref{ex94.14}
we
have seen that
\begin{equation} \label{e94.26}
\Psi(S_l S_m^*) = S_l S_m^* \Lambda_{|m| - |l|}
\end{equation}
(in this and the following equalities we consider $L(l^2(\sZ^+,
l^2(\sZ^+)))$ as an $L(l^2(\sZ^+))$-module in an obvious way).
Thus,
\begin{eqnarray*}
\lefteqn{\Psi(S_l S_m^*) \Psi(S_n S_r^*) - \Psi(S_l S_m^* S_n
S_r^*)} \\
&& = S_l S_m^* S_n S_r^* \left(\Lambda_{|m| - |l|} \Lambda_{|r| -
|n|} - \Lambda_{|r| - |n| + |m| - |r|} \right).
\end{eqnarray*}
For arbitrary integers $a, \, b$ one has
\[
\Lambda_a \Lambda_b - \Lambda_{a+b} = \left\{
\begin{array}{lll}
0 & \mbox{if} & a \le 0, \\
- \Pi_a \Lambda_{a+b} & \mbox{if} & a > 0
\end{array} \right.
\]
with $\Pi_a$ defined as in (\ref{e94.3}). Hence,
\begin{eqnarray*}
\lefteqn{\Psi(S_l S_m^*) \, \Psi(S_n S_r^*) - \Psi(S_l S_m^* S_n
S_r^*)} \\
&& = \left\{
\begin{array}{lll}
- S_l S_m^* S_n S_r^* \, \Pi_{|m| - |l|} \, \Lambda_{|r| - |n| +
|m| - |r|}
  & \mbox{if} & |m| - |l| > 0 \\
0 & \mbox{if} & |m| - |l| \le 0
\end{array} \right. \\
&& = \left\{
\begin{array}{lll}
- \Pi_{|m| - |l|} \, \Psi(S_l S_m^* S_n S_r^*)
  & \mbox{if} & |m| - |l| > 0 \\
0 & \mbox{if} & |m| - |l| \le 0.
\end{array} \right.
\end{eqnarray*}
Since $\Pi_n \in \cC_N$ for every $n \ge 1$ by (\ref{e94.4}), the
inclusion (\ref{e94.24}) follows. \hfill \qed
\section{The lifting homomorphism} \label{s6}
\subsection{The algebra $(e-p_1) \cS_N (e-p_1)$} \label{ss61}
In what follows is will be convenient to compare and to operate
with multi-indices. Given multi-indices $i = (i_1, \, \ldots, \,
i_k) \in \Omega^k$ and $j = (j_1, \, \ldots, \, j_l) \in \Omega^l$
we define their sum as the multi-index
\[
i + j := (i_1, \, \ldots, \, i_k, \, j_1, \, \ldots, \, j_l) \in
\Omega^{k+l}.
\]
Differences of multi-indices $i, \, k$ can be defined only if one
of the multi-indices is a part of the other. Since addition of
multi-indices is not commutative, we consider differences from the
left and from the right. More precisely, we write $i \prec k$ if
there is a multi-index $j$ such that $i+j = k$, and we write $k
\succ j$ if there is a multi-index $i$ with $i+j = k$. The
multi-indices $j$ and $i$ are uniquely determined, and we denote
them by
\[
j := (-i) + k \quad \mbox{and} \quad i := k - j,
\]
respectively. Note that it follows from (\ref{e92.6}) that the
product $s_i^* s_j$ is not zero only if $i \prec j$ or $j \prec
i$. In the first case on gets
\[
s_i^* s_j = s_i^* s_i s_{(-i) + j} = p_{|i|} s_{(-i) + j} =
s_{(-i) + j} p_{|j|},
\]
whereas in the second case
\[
s_i^* s_j = (s_j s_{(-j) + i})^* s_j = s_{(-j) + i}^* s_j^* s_j =
s_{(-j) + i}^* p_{|j|} = p_{|i|} s_{(-j) + i}^*.
\]
\begin{lemma} \label{l95.1}
Let $i, \, j, \, k, \, l$ be multi-indices (not necessarily of the
same length). Then the product
\[
(e - p_1) s_i s_j^* s_k s_l^*
\]
can be written as
\[
(e - p_1) s_r s_t^*
\]
with multi-indices $r$ and $t$ such that $|r| \ge |t|$, or this
product is zero.
\end{lemma}
{\bf Proof.} The product $(e - p_1) s_i s_j^* s_k s_l^*$ is zero,
or one has $j \prec k$ or $k \prec j$. By Lemma \ref{l92.21},
\[
\begin{array}{llll}
(e - p_1) s_i s_j^* s_k s_l^* & = &
\left\{
\begin{array}{l}
(e-p_1) s_i s_j^* s_j s_{(-j)+k} s_l^*   \\
(e-p_1) s_i s_{(-k)+j}^* s_k^* s_k s_l^*
\end{array}
\right. &
\begin{array}{ll}
\mbox{if} & j \prec k \\
\mbox{if} & k \prec j
\end{array} \\
& = &
\left\{
\begin{array}{l}
(e-p_1) s_i p_{|j|} s_{(-j)+k} s_l^*    \\
(e-p_1) s_i s_{(-k)+j}^* p_{|k|} s_l^*
\end{array}
\right. &
\begin{array}{ll}
\mbox{if} & j \prec k \\
\mbox{if} & k \prec j
\end{array} \\
& = &
\left\{
\begin{array}{l}
s_i (p_{|i|} - p_{|i|+1}) p_{|j|} s_{(-j)+k} s_l^*   \\
s_i (p_{|i|} - p_{|i|+1}) p_{|j|} s_{(-k)+j}^*  s_l^*
\end{array}
\right. &
\begin{array}{ll}
\mbox{if} & j \prec k \\
\mbox{if} & k \prec j
\end{array} \\
& = &
\left\{
\begin{array}{l}
0 \\
s_i (p_{|i|} - p_{|i|+1}) s_{(-j)+k} s_l^* \\
s_i (p_{|i|} - p_{|i|+1}) s_{(-k)+j}^*  s_l^*
\end{array}
\right. &
\begin{array}{ll}
\mbox{if} & |j| > |i| \\
\mbox{if} & j \prec k, \, |j| \le |i| \\
\mbox{if} & k \prec j, \, |j| \le |i|
\end{array} \\
& = &
\left\{
\begin{array}{l}
0 \\
(e-p_1) s_i s_{(-j)+k} s_l^*  \\
(e-p_1) s_i s_{(-k)+j}^*  s_l^*
\end{array}
\right. &
\begin{array}{ll}
\mbox{if} & |j| > |i| \\
\mbox{if} & j \prec k, \, |j| \le |i| \\
\mbox{if} & k \prec j, \, |j| \le |i|
\end{array} \\
& = &
\left\{
\begin{array}{l}
0 \\
(e-p_1) s_{i + ((-j)+k)} s_l^* \\
(e-p_1) s_i s_{l + ((-k)+j)}^*
\end{array}
\right. &
\begin{array}{ll}
\mbox{if} & |j| > |i| \\
\mbox{if} & j \prec k, \, |j| \le |i| \\
\mbox{if} & k \prec j, \, |j| \le |i|.
\end{array}
\end{array}
\]
It remains to show that, whenever the product $(e - p_1) s_r
s_t^*$ is not zero, then $|r| \ge |t|$. Assume that $|r| < |t|$.
Then $t$ can be written as $t_1 + t_2$ with $|t_1| > 0$ and $|t_2|
= |r|$. From Lemma \ref{l92.21} $(d)$ we conclude that
\[
(e - p_1) s_r s_t^* = (e - p_1) s_r s_{t_2}^* s_{t_1}^*
= s_r (p_{|r|} - p_{|r| + 1}) s_{t_2}^* s_{t_1}^*
= s_r s_{t_2}^* (e - p_1) s_{t_1}^*,
\]
and $(e - p_1) s_{t_1}^* = 0$ by Lemma \ref{l92.21} $(a)$. \hfill
\qed \\[3mm]
A repeated application of Lemma \ref{l95.1} yields the following.
\begin{coro} \label{c95.2}
Let $i, \, j, \, k, \, l, \, \ldots, \, m, \, n$ be multi-indices
(not necessarily of the same length). Then the product
\[
(e - p_1) s_i s_j^* s_k s_l^* \ldots s_m s_n^* (e - p_1)
\]
can be written as
\[
(e - p_1) s_r s_t^* (e - p_1)
\]
with multi-indices $r$ and $t$ of the same length, or this product
is zero.
\end{coro}
\begin{coro} \label{c95.3}
Let $a \in \cS_N$. Then $(e - p_1) a (e - p_1)$ can be
approximated as closely as desired by linear combinations of
elements of the form
\[
(e - p_1) s_r s_t^* (e - p_1)
\]
with multi-indices $r$ and $t$ of the same length.
\end{coro}
Let $\cS_N^{par}$ stand for the smallest closed subalgebra of
$\cS_N$ which contains all products $s_i s_j^*$ with multi-indices
$i, \, j$ of the same length. Again we allow multi-indices of
length zero, for which we set $s_{\emptyset} := e$.
\begin{lemma} \label{l95.3a}
$\cS_N^{par} = \clos \mbox{\rm span} \:  \{ s_i s_j^* : |i| = |j|
\}$.
\end{lemma}
{\bf Proof.} Let $i, \, j, \, k, \, l$ be multi-indices with $|i|
= |j|$ and $|k| = |l|$. We have to show that the product $(s_i
s_j^*) \, (s_k s_l^*)$ can be written as $s_r s_t^*$ with
multi-indices $r, \, t$ of the same length. This product is zero
if not $j \prec k$ or $k \prec j$. Let, for instance, $j \prec
k$. Then
\begin{eqnarray*}
s_i s_j^* s_k s_l^* & = & s_i s_j^* s_j s_{(-j) + k} s_l^* \\
& = & s_i p_{|i|} s_{(-j) + k} s_l^* \\
& = & s_i s_{(-j) + k} s_l^* \; = \; s_{i + ((-j) + k)} s_l^*
\end{eqnarray*}
where $|i + ((-j) + k)| = |i| + |k| - |j| = |l|$. \hfill \qed
\\[3mm]
One can now state the assertion of Corollary \ref{c95.3} as
follows:
\begin{equation} \label{e95.3b}
\mbox{If} \; a \in \cS_N \quad \mbox{then} \quad (e-p_1) a (e-p_1)
\in (e-p_1) \, \cS_N^{par} \, (e-p_1).
\end{equation}
The mapping
\begin{equation} \label{e95.3c}
\cS_N \to \cS_N^{par}, \quad a \mapsto (e-p_1) a (e-p_1)
\end{equation}
is an expectation which is related with the expectation $\Phi_0 :
\cO_N \to \cO_N^{par}$ as follows.
\begin{prop} \label{p95.3d}
Let $A \in \cO_N$, $a := (P_{N^n} A P_{N^n})_{n \ge 0} +
\cG_\eta$, and $i, \, j \in \sZ^+$. Then
\begin{eqnarray} \label{e95.3e}
\lefteqn{(e-p_1) s_0^i a (s_0^*)^j (e-p_1)} \nonumber \\
&& = (e-p_1) \, \left( (P_{N^n} \Phi_0(S_0^i A (S_0^*)^j)
P_{N^n})_{n \ge 0} + \cG_\eta \right) \, (e-p_1).
\end{eqnarray}
\end{prop}
{\bf Proof.} Because of
\begin{eqnarray*}
s_0^i a (s_0^*)^j & = & (P_{N^n} S_0^i P_{N^n}) \, (P_{N^n} A
P_{N^n}) \, (P_{N^n} (S_0^*)^j P_{N^n}) + \cG_\eta\\
& = & (P_{N^n} S_0^i P_{N^n} A P_{N^n} (S_0^*)^j P_{N^n}) +
\cG_\eta\\
& = & (P_{N^n} S_0^i A (S_0^*)^j P_{N^n}) + \cG_\eta,
\end{eqnarray*}
it is sufficient to prove the assertion for $i = j = 0$. Thus, we
have to show that
\begin{equation} \label{e95.3f}
(e-p_1) a (e-p_1) = (e-p_1) \, \left( (P_{N^n} \Phi_0(A)
P_{N^n})_{n \ge 0} + \cG_\eta \right) \, (e-p_1).
\end{equation}
Both sides of this equality depend linearly and continuously on
$A$. It is thus sufficient to verify (\ref{e95.3f}) for $A = S_k
S_l^*$ with multi-indices $k, l$ of arbitrary length. If $k$ and
$l$ are of the same length, then $A \in \cO_N^{par}$, whence
$\Phi_0(A) = A$, whereas otherwise $\Phi_0(A) = 0$. Thus, the
right-hand side of (\ref{e95.3f}) is equal to $(e-p_1) a (e-p_1)$
if $|k| = |l|$ and zero otherwise. Since
\[
P_{N^n} S_k S_l^* P_{N^n} = P_{N^n} S_k P_{N^n} \cdot P_{N^n}
S_l^* P_{N^n},
\]
the left-hand side of (\ref{e95.3f}) is also zero whenever $|k|
\neq |l|$. \hfill \qed
\subsection{Definition of the lifting homomorphism} \label{ss62}
The desired lifting homomorphism will be defined explicitly by
means of strong limits which involve the following reflection
operators. For every positive integer $n$, let
\[
R_n : l^2(\sZ^+) \to l^2(\sZ^+), \quad (x_k)_{k \ge 0} \mapsto
(x_{n-1}, \, x_{n-2}, \, \ldots, \, x_0, \, 0, \, 0, \, \ldots).
\]
\begin{prop} \label{p95.4}
Let $a \in \cS_N$ and write $(e - p_1) a (e - p_1)$ as $(A_{N^n})
+ \cG_\eta$. Then the strong limit
\begin{equation} \label{e95.5}
\mbox{\rm s-lim}_{n \to \infty} \, R_{N^n} A_{N^n} R_{N^n}
\end{equation}
exists, and the limit is independent of the choice of the
representative of the coset $(e - p_1) a (e - p_1)$.
\end{prop}
{\bf Proof.} Evidently, the limit (\ref{e95.5}) is zero if
$(A_{N^n})$ is a sequence in $\cG_\eta$. This implies the
independence of the choice of the representative.

By Corollary \ref{c95.3}, it is sufficient to prove the existence
of the strong limit (\ref{e95.5}) for sequences $(A_{N^n})$ which
belong to the coset $e - p_1$ or to the coset $s_i s_j^*$ with
multi-indices $i, \, j$ of the same length.

Consider $e - p_1 = (P_{N^n}) - (P_{N^n} S_0^* P_{N^n}) \,
(P_{N^n} S_0 P_{N^n}) + \cG_\eta.$ For every $n \ge 1$, one has
\[
P_{N^n} S_0^* P_{N^n} S_0 P_{N^n} = \diag (1, \, 1, \, \ldots 1,
\, 0, \, 0, \, \ldots, \, 0)
\]
with $N^{n-1}$ ones followed by $N^n - N^{n-1}$ zeros. Hence,
\[
R_{N^n} (P_{N^n} - P_{N^n} S_0^* P_{N^n} S_0 P_{N^n}) R_{N^n}
= \diag (1, \, 1, \, \ldots 1,
\, 0, \, 0, \, \ldots, \, 0)
\]
with $N^n - N^{n-1}$ ones followed by $N^{n-1}$ zeros, which
implies that
\begin{equation} \label{e95.6}
R_{N^n} (P_{N^n} - P_{N^n} S_0^* P_{N^n} S_0 P_{N^n}) R_{N^n}
\to I \quad \mbox{strongly}.
\end{equation}
Next consider the sequence $(P_{N^n} S_i S_j^* P_{N^n})$ with
multi-indices $i, \, j \in \Omega^k$. From Section \ref{ss23} we
infer that
\begin{eqnarray} \label{e95.7}
\lefteqn{S_i S_j^* : \quad (x_r)_{r \ge 0} \mapsto} \nonumber \\
&& ( \underbrace{0, \, \ldots, \, 0,}_{v_{i,k}} \, x_{v_{j,k}}, \,
\underbrace{0, \, \ldots, \, 0,}_{d_k} \, x_{v_{j,k} + N^k}, \,
\underbrace{0, \, \ldots, \, 0,}_{d_k} \, x_{v_{j,k} + 2 N^k}, \,
\ldots).
\end{eqnarray}
In particular, $S_j S_j^*$ is the diagonal projection operator
\begin{equation} \label{e95.8}
S_j S_j^* = \diag (\underbrace{0, \, \ldots, \, 0,}_{v_{j,k}} \,
1, \, \underbrace{0, \, \ldots, \, 0,}_{d_k} \, 1, \,
\underbrace{0, \, \ldots, \, 0,}_{d_k} \, 1, \, \ldots).
\end{equation}
For $n \ge k$, (\ref{e95.7}) implies that
\begin{eqnarray} \label{e95.9}
\lefteqn{P_{N^n} S_i S_j^* P_{N^n} : \quad (x_r)_{r \ge 0} \mapsto}
\nonumber \\
&& ( \underbrace{0, \, \ldots, \, 0,}_{v_{i,k}} \, x_{v_{j,k}}, \,
\underbrace{0, \, \ldots, \, 0,}_{d_k} \, x_{v_{j,k} + N^k}, \,
\ldots, \, x_{v_{j,k} + (N^{n-k} -1) N^k}, \,
\underbrace{0, \, \ldots, \, 0}_{N^k - v_{i,k} -1}).
\end{eqnarray}
Let $V$ denote the shift operator
\[
V : l^2(\sZ^+) \to l^2(\sZ^+), \quad (x_k)_{k \ge 0} \mapsto
(0, \, x_0, \, x_1, \, \ldots).
\]
For every positive integer $n$, set $V_n := V^n$ and $V_{-n} :=
(V^*)^n$, and define $V_0 := I$. Then (\ref{e95.7}) shows that
\[
S_i S_j^* = V_{v_{i,k}} V_{-v_{j,k}} S_j S_j^*,
\]
whence
\[
R_{N^n} P_{N^n} S_i S_j^* P_{N^n} R_{N^n} = R_{N^n} P_{N^n}
V_{v_{i,k}} V_{-v_{j,k}} S_j S_j^* P_{N^n} R_{N^n}.
\]
By (\ref{e95.8}), the operator $S_j S_j^*$ commutes with $P_{N^n}$
which gives
\begin{eqnarray} \label{e95.10}
\lefteqn{R_{N^n} P_{N^n} S_i S_j^* P_{N^n} R_{N^n}} \nonumber \\
&& = R_{N^n} P_{N^n} V_{v_{i,k}} V_{-v_{j,k}} P_{N^n} R_{N^n} \cdot
R_{N^n} P_{N^n} S_j S_j^* P_{N^n} R_{N^n}.
\end{eqnarray}
It is well known and easy to check that
\[
R_{N^n} P_{N^n} V_{v_{i,k}} V_{-v_{j,k}} P_{N^n} R_{N^n} \to
V_{v_{j,k} - v_{i,k}}
\]
strongly as $n \to \infty$. According to (\ref{e95.9}), the second
factor in (\ref{e95.10}) is equal to
\[
R_{N^n} P_{N^n} S_j S_j^* P_{N^n} R_{N^n} =
\diag ( \underbrace{0, \, \ldots, \, 0,}_{N^k - v_{j,k} - 1} \,
1, \, \underbrace{0, \, \ldots, \, 0,}_{d_k} \, 1, \ldots, \, 1, \,
\underbrace{0, \, \ldots, \, 0}_{v_{j,k}}).
\]
Since the number $N^k - v_{j,k} - 1$ of the leading zeros is
independent of $n$, this factor converges strongly to
\[
\diag ( \underbrace{0, \, \ldots, \, 0,}_{N^k - v_{j,k} - 1} \,
1, \, \underbrace{0, \, \ldots, \, 0,}_{d_k} \, 1, \,
\underbrace{0, \, \ldots, \, 0,}_{d_k} \ldots).
\]
Using the notation
\[
\Pi_k := \diag ( 1, \, \underbrace{0, \, \ldots, \, 0,}_{d_k} \, 1, \,
\underbrace{0, \, \ldots, \, 0,}_{d_k} \ldots),
\]
we finally arrive at
\begin{eqnarray} \label{e95.11}
R_{N^n} P_{N^n} S_j S_j^* P_{N^n} R_{N^n} & \to &
V_{v_{j,k} - v_{i,k}} V_{N^k - v_{j,k} - 1} \Pi_k
V_{N^k - v_{j,k} - 1}^* \nonumber \\
& = & V_{N^k - v_{i,k} - 1} \Pi_k V_{N^k - v_{j,k} - 1}^*
\end{eqnarray}
where the latter equality holds since $N^k - v_{j,k} - 1 \ge 0$.
This settles the desired strong convergence. \hfill \qed \\[3mm]
For later use it will prove convenient to write the operator
(\ref{e95.11}) in a different form. Note that
\[
S_i S_j^* = V_{v_{i,k}} V_{-v_{j,k}} S_j S_j^* =
V_{v_{i,k}} V_{-v_{j,k}} V_{v_{j,k}} \Pi_k V_{v_{j,k}}^*
= V_{v_{i,k}} \Pi_k V_{v_{j,k}}^*.
\]
A comparison with (\ref{e95.11}) suggests to introduce the {\em
dual index} $\hat{i}$ of a multi-index $i = (i_1, \, i_2, \,
\ldots, \, i_k)$ by
\[
\hat{i} := (N - 1 - i_1, \, N - 1 - i_2, \, \ldots, \, N - 1 - i_k).
\]
Evidently, $|\hat{i}| = |i| = k$, and one easily checks that
\[
v_{\hat{i}, k} = (N - 1 - i_1) + (N - 1 - i_2) N + \ldots +
(N - 1 - i_k) N^{k-1} = N^k - 1 - v_{i,k}.
\]
Hence,
\[
V_{N^k - v_{i,k} - 1} \Pi_k V_{N^k - v_{j,k} - 1}^* =
V_{v_{\hat{i},k}} \Pi_k V_{v_{\hat{j}, k}}^* =  S_{\hat{i}}
S_{\hat{j}}^*.
\]
For each $i \in \Omega$ we set $S_i^\sharp := S_{\hat{i}} =
S_{N - 1 - i}$. Due to the universal property of Cuntz algebras, the
mapping $^\sharp : S_i \mapsto S_{N - 1 - i}$ can be continued to an
automorphism of $\cO_N$.
\begin{coro} \label{c95.12}
If $i$ and $j$ are multi-indices of the same length then
\begin{equation} \label{e95.13}
R_{N^n} P_{N^n} S_i S_j^* P_{N^n} R_{N^n} \to (S_i S_j^*)^\sharp
\quad \mbox{strongly as} \; n \to \infty.
\end{equation}
\end{coro}
Note that the strong limit (\ref{e95.13}) need not to exist if the
multi-indices are of different length; for example, the sequence
\[
(R_{N^n} P_{N^n} S_0^* P_{N^n} R_{N^n})_{n \ge 0}
\]
for the single isometry $S_0$ does not converge strongly.

We denote the strong limit (\ref{e95.5}) by $W_{00} (a)$ and
consider $W_{00}$ as a mapping from $\cS_N$ into $L(l^2(\sZ^+))$.
More general, for $i, \, j \in \sZ^+$, let $W_{ij} : \cS_N \to
L(l^2(\sZ^+))$ refer to the operator
\[
a \mapsto (S_{N-1}^*)^i \, W_{00} ((e - p_1) s_0^i a (s_0^*)^j (e
- p_1)) S_{N-1}^j
\]
which is consistent with the previous definition. With every
element $a \in \cS_N$, we associate the infinite matrix
\begin{equation} \label{e95.14}
\widetilde{W} (a) := (W_{ij}(a))_{i, \, j \ge 0}
\end{equation}
the entries of which are operators on $L(l^2(\sZ^+))$. We are
going to show that the matrix (\ref{e95.14}) defines a linear
bounded operator on $l^2(\sZ^+, l^2(\sZ^+))$ and that the mapping
$\widetilde{W}$ is the desired injective lifting homomorphism. The
following is the main result of this paper.
\begin{theo} \label{t95.15}
The mapping $\widetilde{W}$ defined by $(\ref{e95.14})$ is a
$^*$-isomorphism from $\cS_N$ onto $\cT_N$ which maps the ideal
$\cJ_N$ onto $\cC_N$.
\end{theo}
\begin{coro} \label{c95.15a}
$\cJ_N$ is the only non-trivial closed ideal of $\cS_N$.
\end{coro}
Indeed, this follows immediately from Theorem \ref{t95.15}
and Corollary \ref{c94.11}.

The proof of Theorem \ref{t95.15} will be given in the
following section. Before, we consider a few examples which
illustrate the action of $\widetilde{W}$. If $a = s_r$ with
$r \in \Omega$, then
\[
(e-p_1) s_0^i s_r (s_0^*)^j (e-p_1) = 0 \quad \mbox{for} \; i+1
\neq j.
\]
If $i+1 = j$ then we get with Corollary \ref{c95.12} that
\begin{eqnarray*}
\lefteqn{(S_{N-1}^*)^i W_{00}((e-p_1) s_0^i s_r (s_0^*)^j (e-p_1))
S_{N-1}^j} \\
&& = (S_{N-1}^*)^i S_{N-1}^i S_{N-1-r} (S_{N-1}^*)^j S_{N-1}^j =
S_{N-1-r} = S_r^\sharp.
\end{eqnarray*}
Thus,
\begin{equation} \label{e95.15a}
\widetilde{W} (s_r) = \Sigma_{N-1-r} \quad \mbox{and, analogously,}
\quad \widetilde{W} (s_r^*) = \Sigma_{N-1-r}^*
\end{equation}
with $\Sigma_i$ defined by (\ref{e94.2}). Let now $a = s_0^*s_0 = p_1$.
If $i=0$, one has for every $j$
\[
(e-p_1) a (s_0^*)^j (e-p_1) = (e-p_1) p_1 (s_0^*)^j (e-p_1) = 0,
\]
and the same result follows if $i > 0$ and $i\neq j$:
\[
(e-p_1) s_0^i s_0^* s_0 (s_0^*)^j (e-p_1) = (e-p_1) s_0^i (s_0^*)^j
(e-p_1) = 0.
\]
Let, finally, $i > 0$ and $i=j$. Then
\[
(e-p_1) s_0^i s_0^* s_0 (s_0^*)^i (e-p_1) = (e-p_1) s_0^i (s_0^*)^i
(e-p_1).
\]
Applying the homomorphism $W_{00}$ to the right-hand side of this
equality and taking into account Corollary \ref{c95.12} we obtain the
operator $S_{N-1}^i (S_{N-1}^*)^i$, whence
\[
\widetilde{W} (p_1) = I - \Pi_1 \quad \mbox{and} \quad \widetilde{W}
(e-p_1) = \Pi_1.
\]
Note that $\widetilde{W} (s_r^*) \widetilde{W} (s_r) = I - \Pi_1 =
\widetilde{W} (p_1) = \widetilde{W} (s_r^* s_r)$ in contrast to the
mapping $\Psi$ which does not act multiplicatively on these elements.
\subsection{Proof of Theorem \ref{t95.15}} \label{ss63}
It is not too hard to verify that the mapping $\widetilde{W}$ acts as a $^*$-homomorphism on a {\em dense subalgebra} of $\cS_N$. That it acts as a $^*$-homomorphism on the whole algebra would easily follow from this fact if one would know that $\widetilde{W}$ is bounded. Conversely, the boundedness of $\widetilde{W}$ comes as a simple consequence of the fact that $\widetilde{W}$ acts as a $^*$-homomorphism on $\cS_N$. Unfortunately, neither the boundedness of $\widetilde{W}$ nor the homomorphy of $\widetilde{W}$ on all of $\cS_N$ could be shown
directly. Rather we have to prove both properties simultaneously by climbing step by step form small substructures of $\cS_N$ to the whole algebra.

We will make use of the fact that the algebra $\cS_N$ splits into
the direct sum
\begin{equation} \label{e95.16}
\{ (P_{N^n} A P_{N^n})_{n \ge 0} + \cG_\eta : A \in \cO_N \} \oplus
\cJ_N.
\end{equation}
\begin{prop} \label{p95.17}
The mapping $\widetilde{W}$ acts as a linear contraction on the first
summand of $(\ref{e95.16})$, and it maps this summand into $\cT_N$.
\end{prop}
{\bf Proof.} Let $A \in \cO_N$ and $a := (P_{N^n} A P_{N^n})_{n \ge 0}
+ \cG_\eta$. From equality (\ref{e95.3e}) we conclude that
\begin{eqnarray*}
\lefteqn{W_{00} ((e-p_1) s_0^i a (s_0^*)^j (e-p_1))} \\
&& = W_{00} \left((e-p_1) \, \left( (P_{N^n} \Phi_0(S_0^i A (S_0^*)^j)
P_{N^n})_{n \ge 0} + \cG_\eta \right) \, (e-p_1) \right) \\
&& = \left( \Phi_0 (S_0^i A (S_0^*)^j) \right)^\sharp
\end{eqnarray*}
Thus,
\begin{eqnarray} \label{e95.18}
\widetilde{W}(a) & = &
\left( (S_{N-1}^*)^i \left( \Phi_0 (S_0^i A (S_0^*)^j) \right)^\sharp
\, S_{N-1}^j \right)_{i,j \ge 0} \nonumber \\
& = & \left( \left( (S_0^*)^i  \, \Phi_0 (S_0^i A (S_0^*)^j) \, S_0^j
\right)^\sharp \, \right)_{i,j \ge 0} \nonumber \\
& = & \left( (S_{N-1}^*)^i \, \Phi_0 (S_0^i A (S_0^*)^j)^\sharp \,
S_{N-1}^j \right)_{i,j \ge 0}.
\end{eqnarray}
We claim that
\[
\Phi_0 (B)^\sharp = \Phi_0(B^\sharp) \quad \mbox{for every} \; B \in
\cO_N.
\]
Since $\cO_N$ is spanned by products $S_l S_m^*$ with multi-indices
$l, \, m$ and since the mappings $\Phi_0$ and $^\sharp$ are continuous,
it is sufficient to check the claim for $B = S_l S_m^*$. If $|l| \neq
|m|$, then both sides of the claimed identity are zero, whereas
\[
\Phi_0 (S_l S_m^*)^\sharp = (S_l S_m^*)^\sharp = S_{\hat{l}} S_{\hat{m}}^*
= \Phi_0 (S_{\hat{l}} S_{\hat{m}}^*) = \Phi_0 ((S_l S_m^*)^\sharp)
\]
if $|l| = |m|$. This proves the claim and shows that (\ref{e95.18}) is
equal to
\[
\left( (S_{N-1}^*)^i \, \Phi_0 (S_{N-1}^i A^\sharp (S_{N-1}^*)^j) \,
S_{N-1}^j \right)_{i,j \ge 0}.
\]
Repeating the arguments from Example \ref{ex94.14} one easily
gets that the operators
\[
(S_r^*)^i \, \Phi_0 (S_r^i A^\sharp (S_r^*)^j) \, S_r^j
\]
are independent of the choice of $r \in \Omega$. Thus, (\ref{e95.18})
further coincides with
\[
\left( (S_0^*)^i \, \Phi_0 (S_0^i A^\sharp (S_0^*)^j) \,
S_0^j \right)_{i,j \ge 0}
\]
whence
\[
\widetilde{W}(a) = \Psi(A^\sharp).
\]
Now we conclude as follows. The mapping
\[
a = (P_{N^n} A P_{N^n}) + \cG_\eta \; \mapsto \; A = \mbox{s-lim}
P_{N^n} A P_{N^n}
\]
is a linear contraction by the Banach-Steinhaus theorem. As already
mentioned, the mapping $A \mapsto A^\sharp$ is a linear contraction
(and even an isometry) on $\cO_N$, and from Theorem \ref{t94.20} we
recall that the mapping $A^\sharp \mapsto \Psi(A^\sharp)$ is a linear
contraction, too, with range in $\cT_N$. \hfill \qed \\[3mm]
Now we consider the second summand in (\ref{e95.16}). Abbreviate
$e - p_n$ to $\pi_n$ and recall the definition (\ref{e94.3}) of
$\Pi_n$.
\begin{prop} \label{p95.19}
The mapping $\widetilde{W}$ is a $^*$-homomorphism from $\pi_n \cJ_N
\pi_n$ into $\Pi_n \cC_N \Pi_n$ for every $n \ge 1$.
\end{prop}
{\bf Proof.} Let $a \in \cJ_N$. First we show that $\widetilde{W}
(\pi_n a \pi_n) \in \Pi_n \cC_N \Pi_n$. For, write $\pi_n = e - p_n
= \sum_{i=0}^{n-1} (p_i - p_{i+1})$ with the convention $p_0 := e$.
Then
\begin{eqnarray} \label{e95.20}
\pi_n a \pi_n & = & \sum_{i,k = 0}^{n-1} (p_i - p_{i+1}) a (p_k -
p_{k+1}) \nonumber \\
& = & \sum_{i,k = 0}^{n-1} (s_0^*)^i (e-p_1) s_0^i a (s_0^*)^k (e-p_1)
s_0^k \nonumber \\
& = & \sum_{i,k = 0}^{n-1} (s_0^*)^i (e-p_1) d_{ik} (e-p_1) s_0^k
\end{eqnarray}
where the $d_{ik}$ can be found in $\cS_N^{par}$ by Corollary
\ref{c95.3} and (\ref{e95.3b}). Since $\cS_N^{par}$ is spanned by
the products $s_\mu s_\nu^*$ with multi-indices $\mu$ and $\nu$ of
the same length, the assertion will follow once we have shown that
\begin{equation} \label{e95.21}
\widetilde{W} ((s_0^*)^i (e-p_1) s_\mu s_\nu^* (e-p_1) s_0^k)
\in \Pi_n \cC_N \Pi_n
\end{equation}
whenever $|\mu| = |\nu|$. The right-hand side of (\ref{e95.21})
is the matrix
\[
\left( (S_{N-1}^*)^r W_{00} ( (e-p_1) s_0^r (s_0^*)^i (e-p_1)
s_\mu s_\nu^* (e-p_1) s_0^k (s_0^*)^t (e-p_1) ) S_{N-1}^t \right)
_{r,t \ge 0}.
\]
Only the $ik$th entry of this matrix is non-zero, and this entry
equals
\[
(S_{N-1}^*)^i W_{00} ( (e-p_1) s_0^i (s_0^*)^i (e-p_1)
s_\mu s_\nu^* (e-p_1) s_0^k (s_0^*)^k (e-p_1)) S_{N-1}^k.
\]
By Proposition \ref{p95.4}, this entry further coincides with
\begin{eqnarray*}
\lefteqn{(S_{N-1}^*)^i  \, W_{00} ( (e-p_1) s_0^i (s_0^*)^i
(e-p_1)) \times} \\
&& \times W_{00} ( (e-p_1) s_\mu s_\nu^* (e-p_1)) \,
W_{00} ( (e-p_1) s_0^k (s_0^*)^k (e-p_1)) \, S_{N-1}^k
\end{eqnarray*}
which on its hand is the same as
\[
(S_{N-1}^*)^i S_{N-1}^i (S_{N-1}^*)^i S_{\hat{\mu}} S_{\hat{\nu}}^*
S_{N-1}^k (S_{N-1}^*)^k S_{N-1}^k = (S_{N-1}^*)^i S_{\hat{\mu}}
S_{\hat{\nu}}^* S_{N-1}^k
\]
by Corollary \ref{c95.12}. Thus, the right-hand side of
(\ref{e95.21}) coincides with
\[
(\Pi_i - \Pi_{i-1}) \Sigma_{N-1} \Psi(S_{\hat{\mu}}
S_{\hat{\nu}}^*) \Sigma_{N-1}^* (\Pi_k - \Pi_{k-1})
\]
(with the convention $\Pi_0 := 0$). Clearly, this matrix is in
$\Pi_n \cC_N \Pi_n$, which proves (\ref{e95.21}).

It is now evident that the mapping $\widetilde{W} : \pi_n \cJ_N
\pi_n \to \Pi_n \cC_N \Pi_n$ is linear and symmetric. It remains
to verify that this mapping is multiplicative,
\begin{equation} \label{e95.22}
\widetilde{W} (\pi_n a \pi_n) \, \widetilde{W} (\pi_n b \pi_n)
= \widetilde{W} (\pi_n a \pi_n b \pi_n) \quad \mbox{for} \; a, \, b \in
\cJ_N.
\end{equation}
It is an elementary fact that a linear mapping $W$ between algebras
$\cA$ and $\cB$ which are spanned by subsets $\cL$ and $\cM$,
respectively, is multiplicative if (and only if) $W(l) W(m) = W(lm)$
for each choice of elements $l \in \cL$ and $m \in \cM$. Thus, and by
the first part of this proof, it is sufficient to verify the equality
(\ref{e95.22}) in case
\[
\pi_n a \pi_n = (s_0^*)^i (e-p_1) s_\mu s_\nu^* (e-p_1) s_0^k,
\]
\[
\pi_n b \pi_n = (s_0^*)^l (e-p_1) s_\lambda s_\tau^* (e-p_1) s_0^m
\]
with multi-indices $\mu, \, \nu, \, \lambda$ and $\tau$ such that
$|\mu| = |\nu|$ and $|\lambda| = |\tau|$. From above we infer that
only the $ik$th entry of the matrix $\widetilde{W} (\pi_n a \pi_n)$
is different from zero, and this entry is $(S_{N-1}^*)^i S_{\hat{\mu}}
S_{\hat{\nu}}^* S_{N-1}^k$. Similarly, only the $lm$th entry of
the matrix $\widetilde{W} (\pi_n b \pi_n)$ not zero, and this entry
equals $(S_{N-1}^*)^l S_{\hat{\lambda}} S_{\hat{\tau}}^* S_{N-1}^m$.
Consequently, the matrix $\widetilde{W} (\pi_n a \pi_n) \,
\widetilde{W} (\pi_n b \pi_n)$ is not the zero matrix only if $k=l$.
In this case, this matrix has at most one non-vanishing entry,
namely the $im$th entry, which is
\begin{equation} \label{e95.23}
(S_{N-1}^*)^i S_{\hat{\mu}} S_{\hat{\nu}}^* S_{N-1}^k \cdot
(S_{N-1}^*)^k S_{\hat{\lambda}} S_{\hat{\tau}}^* S_{N-1}^m.
\end{equation}
Now we consider $\widetilde{W} (\pi_n a \pi_n b \pi_n)$. It is
$\pi_n a \pi_n b \pi_n \neq 0$ only if $k = l$, in which case
\[
\pi_n a \pi_n b \pi_n = (s_0^*)^i (e-p_1) s_\mu s_\nu^* (e-p_1)
s_0^k \cdot (s_0^*)^k (e-p_1) s_\lambda s_\tau^* (e-p_1) s_0^m.
\]
Further, the product
\[
(e-p_1) s_0^r \pi_n a \pi_n b \pi_n (s_0^*)^t (e-p_1)
\]
is not zero only if $r=i$ and $m=t$. Consequently, the matrix
$\widetilde{W} (\pi_n a \pi_n b \pi_n)$ is not zero only if
$k=l$. In this case, only the $im$th entry of this matrix does
not vanish. This entry is
\[
(S_{N-1}^*)^i S_{N-1}^i (S_{N-1}^*)^i S_{\hat{\mu}}
S_{\hat{\nu}}^* S_{N-1}^k \cdot (S_{N-1}^*)^k S_{\hat{\lambda}}
S_{\hat{\tau}}^* S_{N-1}^m (S_{N-1}^*)^m S_{N-1}^m,
\]
which coincides with (\ref{e95.23}). \hfill \qed
\begin{coro} \label{c95.24}
The mapping $\widetilde{W}$ is a linear contraction from $\cJ_N$
into $\cC_N$.
\end{coro}
{\bf Proof.} Let $j \in \cJ_N$. From Lemma \ref{l92.23} we infer
that
\begin{equation} \label{e95.25}
\lim_{n \to \infty} \|j - \pi_n j \pi_n\| = 0 \quad \mbox{for each}
\; j \in \cJ_N.
\end{equation}
Hence, $(\pi_n j \pi_n)_{n \ge 1}$ is a Cauchy sequence. From
Proposition \ref{p95.19} we further infer that $\widetilde{W} :
\pi_n \cJ_N \pi_n \to \Pi_n \cC_N \Pi_n$ is a $^*$-homomorphism,
hence contractive:
\begin{equation} \label{e95.26}
\|\widetilde{W} (\pi_n j \pi_n)\| \le \|\pi_n j \pi_n\| \quad
\mbox{for all} \, n \ge 1.
\end{equation}
Since $\pi_n j \pi_n \in \pi_m \cJ_N \pi_m$ for $m \ge n$, we
conclude that
\[
\|\widetilde{W} (\pi_n j \pi_n) - \widetilde{W} (\pi_m j \pi_m)\|
\le \|\pi_n j \pi_n - \pi_m j \pi_m\|
\]
whenever $m \ge n$. Hence, $(\widetilde{W} (\pi_n j \pi_n))
_{n \ge 1}$ is a Cauchy sequence. Let $J$ denote its limit. Since
all entries of the matrix mapping $\widetilde{W}$ are continuous,
we conclude from
\[
\|\widetilde{W} (\pi_n j \pi_n) - J\| \to 0
\]
that $J = \widetilde{W} (j)$. Now it is clear that $\widetilde{W}
(j) \in \cC_N$, and passing to the limit as $n \to \infty$ in
(\ref{e95.26}) yields $\|\widetilde{W} (j)\| \le \|j\|$ for every
$j \in cJ_N$. \hfill \qed
\begin{coro} \label{c95.27}
The mapping $\widetilde{W} : \cJ_N \to \cC_N$ is a $^*$-homomorphism.
\end{coro}
{\bf Proof.} We have to show that $\widetilde{W}$ is a multiplicative
mapping on $\cJ_N$. Let $j_1, \, j_2 \in \cJ_N$. By Lemma
\ref{l92.23}, $j_1j_2 = \lim \pi_n j_1 \pi_n j_2 \pi_n$, and since
$\widetilde{W}$ is continuous on $\cJ_N$,
\[
\widetilde{W} (j_1j_2) = \lim \widetilde{W} (\pi_n j_1 \pi_n
j_2 \pi_n).
\]
Since $\widetilde{W}$ is multiplicative on $\pi_n \cJ_N \pi_n$
by Proposition \ref{p95.19}, this implies
\[
\widetilde{W} (j_1j_2) = \lim \widetilde{W} (\pi_n j_1 \pi_n)
\, \widetilde{W} (\pi_n j_2 \pi_n) = \widetilde{W} (j_1)
\widetilde{W}(j_2),
\]
whence the assertion. \hfill \qed
\begin{coro} \label{c95.28}
The mapping $\widetilde{W}$ is bounded on all of $\cS_N$.
\end{coro}
{\bf Proof.} Let $(A_n) + \cG_\eta \in \cS_N$.  In accordance
with (\ref{e95.16}), we write this coset as
\begin{equation} \label{e95.29}
(A_n) + \cG_\eta = \left( (P_{N^n} A P_{N^n}) + \cG_\eta
\right) + \left( (J_n) + \cG_\eta \right) =: a + j
\end{equation}
with $A := \mbox{s-lim} A_n P_{N^n}$ and $j = (J_n) + \cG_\eta
\in \cJ_N$. Then
\[
\|a\| = \| (P_{N^n} A P_{N^n}) + \cG_\eta \| \le \|A\| \le
\|(A_n) + \cG_\eta\|,
\]
i.e., the first summand in (\ref{e95.29}) depends continuously
on $(A_n) + \cG_\eta$. From $\|a\| \le \|a+j\|$ we obtain
$\|j\| \le \|a+j\| + \|a\| \le 2 \, \|a+j\|$, whence
\[
\|\widetilde{W}(a+j)\| \le \|\widetilde{W}(a)\| + \|\widetilde{W}(j)\|
\le \|a\| + \|j\| \le 3 \, \|a+j\|
\]
due to Proposition \ref{p95.17} and Corollary \ref{c95.24}. \hfill
\qed
\begin{prop} \label{p95.30}
The mapping $\widetilde{W}$ is a $^*$-homomorphism from $\cS_N$
into $\cT_N$.
\end{prop}
{\bf Proof.} It is sufficient to verify that $\widetilde{W}$ is
multiplicative on $\cS_N$. We start with showing a partial
multiplicativity result,
\begin{equation} \label{e95.31}
\widetilde{W} (a \pi_k) = \widetilde{W} (a) \Pi_k \quad \mbox{for
each} \; a \in \cS_N.
\end{equation}
Indeed, the matrix representation of $\widetilde{W} (a \pi_k)$ is
\[
\left( (S_{N-1}^*)^i W_{00}( (e-p_1) s_0^i a \underbrace{\pi_k
(s_0^*)^j (e-p_1)}) S_{N-1}^j \right)_{i,j \ge 0},
\]
and the underbraced expression equals
\begin{eqnarray*}
\pi_k (s_0^*)^j (e-p_1) & = & (e-p_k) (s_0^*)^j (e-p_1)\\
& = & (e-p_k) (p_j - p_{j+1}) (s_0^*)^j \\
& = & \left\{
\begin{array}{lll}
0 & \mbox{if} & k \le j \\
(p_j - p_{j+1} - p_k + p_k p_{j+1}) (s_0^*)^j & \mbox{if} & k > j
\end{array} \right. \\
& = & \left\{
\begin{array}{lll}
0 & \mbox{if} & k \le j \\
(p_j - p_{j+1}) (s_0^*)^j & \mbox{if} & k > j
\end{array} \right. \\
& = & \left\{
\begin{array}{lll}
0 & \mbox{if} & k \le j \\
(s_0^*)^j (e-p_1) & \mbox{if} & k > j
\end{array} \right.
\end{eqnarray*}
whence the assertion (\ref{e95.31}). For the proof of the
general assertion, let $a, \, b \in \cS_N$. Since $\pi_k \in
\cJ_N$ for every $k$ by Proposition \ref{p92.22} and $\widetilde{W}$
is multiplicative on $\cJ_N$, we get
\[
\widetilde{W} ( a \pi_m b \pi_k) = \widetilde{W} ( a \pi_m) \,
\widetilde{W} ( b \pi_k)
\]
for all $k, \, m \ge 0$. By (\ref{e95.31}),
\[
\widetilde{W} (a \pi_m b \pi_k) = \widetilde{W} (a) \Pi_m \,
\widetilde{W} (b) \Pi_k.
\]
For $m \to \infty$ we have $a \pi_m b \pi_k \to ab \pi_k$ by
Lemma \ref{l92.23} (recall that $\pi_k \in \cJ_N$) and
$\Pi_m \to I$ strongly. Thus, due to the continuity of
$\widetilde{W}$,
\[
\widetilde{W} (a b \pi_k) = \widetilde{W} (a) \,
\widetilde{W} (b) \Pi_k.
\]
Invoking (\ref{e95.31}) again and letting $k$ tend to infinity,
we arrive at the assertion. \hfill \qed \\[3mm]
We prepare the proof of the next proposition by a simple lemma.
\begin{lemma} \label{l95.32}
Every coset in $(e-p_1) \cS_N (e-p_1)$ can be written in the form
$(e-p_1) c (e-p_1)$ with $c = (P_{N^n} C P_{N^n})_{n \ge 0} +
\cG_\eta$ with $C \in \cO_N^{par}$.
\end{lemma}
{\bf Proof.} The assertion holds for cosets of the form
$(e - p_1) s_r s_t^* (e - p_1)$ with multi-indices $r, \, t$
of the same length. Indeed, it follows immediately from
the identities (\ref{e92.4}) that
\[
(e - p_1) s_r s_t^* (e - p_1) = (e - p_1) \, \left(
(P_{N^n} S_r S_t^* P_{N^n})_{n \ge 0} + \cG_\eta \right) \,
(e - p_1).
\]
Then, by Corollary \ref{c95.3}, the assertion holds for all
cosets in a dense subalgebra of $(e-p_1) \cS_N (e-p_1)$. Let
now $a$ be an arbitrary element of $\cS_N$. Let $((e-p_1)
a_n (e-p_1))_{n \ge 1}$ be a sequence in this dense subalgebra
which converges to $(e-p_1) a (e-p_1)$ in the norm. As we have
just checked, there are operators $A_n \in \cO_N^{par}$ such
that
\[
(e-p_1) a_n (e-p_1) = (P_{N^n} A_n P_{N^n})_{n \ge 0} +
\cG_\eta.
\]
Applying the homomorphism $W_{00}$ to both sides of this
equality we obtain
\[
W_{00} ((e-p_1) a_n (e-p_1)) = A_n^\sharp
\]
for every $n$. Since the sequence $((e-p_1) a_n (e-p_1))
_{n \ge 1}$ converges, $(A_n^\sharp)$ is a Cauchy sequence.
But then $(A_n)$ is a Cauchy sequence; hence, convergent, too.
Let $A$ denote its limit. Then $A \in \cO_N^{par}$ and, clearly,
$(e-p_1) a (e-p_1) = (P_{N^n} A P_{N^n})_{n \ge 0} + \cG_\eta$.
\hfill \qed
\begin{prop} \label{p95.33}
The mapping $\widetilde{W}$ is injective on $\cJ_N$.
\end{prop}
{\bf Proof.} The assertion will follow once we have shown that
\begin{equation} \label{e95.34}
\widetilde{W} \; \mbox{is injective on} \; \pi_n \cJ_N \pi_n \;
\mbox{for every} \; n \ge 1.
\end{equation}
Indeed, let (\ref{e95.34}) be satisfied, and let $j \in \cJ_N$
be an element with $\widetilde{W} (j) = 0$. Then $\Pi_n
\widetilde{W} (j) \Pi_n =0$ for every $n$. By (\ref{e95.31}),
this implies $\widetilde{W} (\pi_n j \pi_n) = 0$ for every $n$.
From (\ref{e95.34}) we infer that $\pi_n j \pi_n = 0$ for every
$n$. Passage to the limit $n \to \infty$ yields $j=0$, which
implies the desired injectivity.

Further, (\ref{e95.34}) will follow once we have shown that
\begin{equation} \label{e95.35}
\widetilde{W} \; \mbox{is injective on} \; \pi_1 \cJ_N \pi_1
= (e-p_1) \cJ_N (e-p_1).
\end{equation}
Indeed, let (\ref{e95.35}) be satisfied, and let $j \in \cJ_N$
be such that
\[
\widetilde{W} (\pi_n j \pi_n) = \Pi_n \widetilde{W}
(j) \Pi_n = 0 \quad \mbox{for some} \; n \ge 0.
\]
Then
\[
(S_{N-1}^*)^i \, W_{00} ( (e-p_1) s_0^i j (s_0^*)^k (e-p_1)) \,
S_{N-1}^k = 0
\]
for all $i, \, k \le n-1$. Multiplication by $S_{N-1}^i$ from
the left and by $(S_{N-1}^*)^k$ from the right-hand side yields
\begin{equation} \label{e95.36}
S_{N-1}^i (S_{N-1}^*)^i \, W_{00} ( (e-p_1) s_0^i j (s_0^*)^k
(e-p_1)) \, S_{N-1}^k (S_{N-1}^*)^k = 0
\end{equation}
for $i, \, k \le n-1$. Because of
\[
S_{N-1}^i (S_{N-1}^*)^i = W_{00} ( (e-p_1) s_0^i (s_0^*)^i (e-p_1))
\]
and by the multiplicativity of $W_{00}$ on $(e-p_1) \cJ_N (e-p_1)$,
we conclude from (\ref{e95.36}) that
\[
W_{00} ( (e-p_1) s_0^i (s_0^*)^i (e-p_1) s_0^i j (s_0^*)^k
(e-p_1) s_0^k (s_0^*)^k (e-p_1)) = 0.
\]
Since
\[
(e-p_1) s_0^i (s_0^*)^i (e-p_1) = (e-p_1) s_0^i (s_0^*)^i,
\]
this finally implies that
\[
W_{00} ( (e-p_1) s_0^i j (s_0^*)^k (e-p_1)) = 0 \quad \mbox{if} \;
i, \, k \le n-1.
\]
Via assumption (\ref{e95.35}), this gives
\begin{equation} \label{e95.37}
(e-p_1) s_0^i j (s_0^*)^k (e-p_1) = 0 \quad \mbox{if} \;
i, \, k \le n-1
\end{equation}
(note that
\[
\widetilde{W} ((e-p_1) r (e-p_1)) =
\pmatrix{W_{00} ((e-p_1) r (e-p_1)) & 0 & \ldots \cr
0 & 0 & \ldots \cr
\vdots & \vdots & \ddots}
\]
for every element $r$ of $\cJ_N$, which implies that $W_{00}$ is
injective whenever (\ref{e95.35}) holds). Multiplying
(\ref{e95.36}) from the left by $(s_0^*)^i$ and by $s_0^k$ from
the right-hand side and taking into account that $(s_0^*)^i
(e-p_1) s_0^i = p_i - p_{i+1}$ we obtain
\[
(p_i - p_{i+1}) j (p_k - p_{k+1}) = 0 \quad \mbox{if} \;
i, \, k \le n-1.
\]
Summation over $0 \le i, \, k \le n-1$ gives $\pi_n j \pi_n = 0$,
whence the injectivity of $\widetilde{W}$ on $\pi_n \cJ_N \pi_n$.

It remains to prove (\ref{e95.35}). Let $j \in \cJ_N$ and
$\widetilde{W} ((e-p_1) j (e-p_1)) = 0$. Employing Lemma \ref{l95.32}, we can write $(e-p_1) j (e-p_1)$ as
\[
(e-p_1) \left( (P_{N^n} C P_{N^n})_{n \ge 0} + \cG_\eta \right)
(e-p_1) \quad \mbox{with} \; C \in \cO_N^{par}.
\]
Consequently,
\[
0 = W_{00} ((e-p_1) j (e-p_1)) = \Phi_0(C^\sharp) = C^\sharp
\]
since $\Phi$ is an expectation from $\cO_N$ onto $\cO_N^{par}$ and
$C$ belongs to the latter subalgebra. Thus, $C=0$, which implies
that $(e-p_1) j (e-p_1)$. The injectivity of $\widetilde{W}$ on
$\pi_1 \cJ_N \pi_1$ follows. \hfill \qed \\[3mm]
Now we can finish the proof of Theorem \ref{t95.15} as follows.
The mapping $\widetilde{W}$ is an injective $^*$-homomorphism on
$\cJ_N$ as we have just seen. Hence, by Corollary \ref{c93.8},
$\widetilde{W}$ is an injective $^*$-homomorphism on $\cS_N$. The
range of this homomorphism contains the generating operators
$\Sigma_k$ of $\cT_N$; thus, $\widetilde{W}$ maps $\cS_N$ onto
$\cT_N$. Since $\widetilde{W}$ maps the generating element $e-p_1$
of the ideal $\cJ_N$ to the generating element $\Pi_1$ of $\cC_N$,
it is further clear that $\widetilde{W}$ maps $\cJ_N$ onto
$\cC_N$. \hfill \qed
\subsection{Some consequences of Theorem \ref{t95.15}} \label{ss64}
\paragraph{Stability.} The assertion of Theorem \ref{t95.15} is
equivalent to the following stability criterion.
\begin{theo} \label{t95.38}
A sequence $\bA = (A_n)$ in $\cS_\eta(\cO_N)$ is stable if and
only if the operator $\widetilde{W} (\bA + \cG_\eta)$ is invertible.
\end{theo}
Specifying this result to finite sections sequences for operators
in the Cuntz algebra yields
\begin{coro} \label{c95.39}
Let $A \in \cO_N$. Then the sequence $(P_{N^n} A P_{N^n})_{n \ge 0}$ is stable if and only if the (stratified) Toeplitz operator $\Psi(A^\sharp) = T(f_{A^\sharp}) \in L(l^2(\sZ^+, l^2(\sZ^+)))$ is invertible.
\end{coro}
{\bf Proof.} Set $a := (P_{N^n} A P_{N^n})_{n \ge 0} + \cG_\eta$. We have to show that $\widetilde{W} (a) = \Psi(A^\sharp)$. If $A$ is a product $S_l S_m^*$ with multi-indices $l$ and $m$, then this equality follows from (\ref{e94.26}) and (\ref{e95.15a}) by
noting that the latter implies
\[
\widetilde{W} (s_l s_m^*) = S_{\hat{l}} S_{\hat{m}}^*
\Lambda_{|\hat{m}| - |\hat{l}|}
\]
due to the homomorphy of $\widetilde{W}$ and that
\[
S_{\hat{l}} S_{\hat{m}}^* \Lambda_{|\hat{m}| - |\hat{l}|}
= (S_l S_m^*)^\sharp \Lambda_{|m| - |l|}
\]
with $\Lambda$ defined by (\ref{e94.25}). The general case
follows from this partial result since the products $S_l S_m^*$
span a dense subalgebra of $\cO_N$. \hfill \qed
\paragraph{Fractality.} The following is certainly the most
important consequence of Theorem \ref{t95.15}. It can also
serve as a perfect illustration to Theorem \ref{t91.11}. The
proof will follow directly from the special form of the
homomorphism $\widetilde{W}$.
\begin{coro} \label{c95.40}
The algebra $\cS_\eta(\cO_N)$ is fractal.
\end{coro}
{\bf Proof.} Recall that the entries of the matrix operator
$\widetilde{W} ((A_n) + \cG_\eta)$ are defined by strong limits.
Consequently, if only an (infinite) subsequence of $(A_n)$ is
known, one can nevertheless determine the operator $\widetilde{W}
((A_n) + \cG_\eta) \in \cT_N$. Since $\widetilde{W} : \cS_N
\to \cT_N$ is an isomorphism one can, thus, reconstruct the
coset of $(A_n)$ modulo $\cG_\eta$ from each subsequence of
$(A_n)$.  \hfill \qed
\paragraph{Spectral approximation.} As already mentioned,
sequences in fractal algebras are distinguished by their
excellent convergence properties. To mention only a few of
them, let $\sigma (a)$ denote the spectrum of an element
$a$ of a $C^*$-algebra with identity element $e$, write
$\sigma_2(a)$ for the set of the singular values of $a$, i.e.,
$\sigma_2(a)$ is the set of all non-negative square roots
of elements in the spectrum of $a^*a$ and finally, for $\varepsilon
> 0$, let $\sigma^{(\varepsilon)} (a)$ refer to the
$\varepsilon$-pseudospectrum of $a$, i.e. to the set of all
$\lambda \in \sC$ for which $a - \lambda e$ is not invertible or
$\|(a - \lambda e)^{-1}\| \ge 1/\varepsilon$. Let further
\[
d_H(M, \, N) := \max \, \{ \max_{m \in M} \min_{n \in N} |m-n|, \,
\max_{n \in N} \min_{m \in M} |m-n| \}
\]
denote the Hausdorff distance between the non-empty compact subsets
$M$ and $N$ of the complex plane.
\begin{theo} \label{t95.41}
Let $(A_n)$ be a sequence in $\cS_\eta(\cO_N)$ and set $a :=
(A_n) + \cG_\eta$. Then the following set-sequences converge with
respect to the Hausdorff distance as $n \to \infty \!:$ \\[1mm]
$(a)$ $\sigma(A_n) \to \sigma (\widetilde{W}(a))$ if
$a$ is self-adjoint; \\[1mm]
$(b)$ $\sigma_2 (A_n) \to \sigma_2 (\widetilde{W}(a))$;
\\[1mm]
$(c)$ $\sigma^{(\varepsilon)} (A_n) \to \sigma^{(\varepsilon)}
(\widetilde{W}(a))$.
\end{theo}
The proof follows immediately from the stability criterion in
Theorem \ref{t95.15} above and from Theorems 3.20, 3.23 and 3.33
in \cite{HRS2}. Let us emphasize that in general one cannot
remove the assumption $a = a^*$ in assertion $(a)$, whereas $(c)$
holds without any assumption. This observation is only one
reason for the present increasing interest in pseudospectra.
For detailed presentations of pseudospectra and their applications
as well as of other spectral quantities see the monographs
\cite{BGr5,BSi2,HRS2,TrE1} and the references therein.
\paragraph{Compactness and Fredholm properties.} Recall the
definition of the algebra $\cF$ of all bounded sequences of
matrices and of its ideal $\cK$ ideal of the compact sequences
from Section \ref{ss33} and let $\cF_\eta$ and $\cK_\eta$ denote
the corresponding restricted algebras.
\begin{prop} \label{p95.42}
The only compact sequences in $\cS_\eta (\cO_N)$ are the
sequences in $\cG_\eta$.
\end{prop}
{\bf Proof.} By Corollary \ref{c95.15a}, $\cJ_N$ is the only
non-trivial closed ideal of $\cS_N$. Thus, the intersection
$\cS_N \cap (\cK_\eta/\cG_\eta)$ is either $\cS_N, \, \cJ_N,$
or $\{0\}$. Since already $\cJ_N$ contains cosets of non-compact
sequences (e.g., the coset $e - p_1$), the assertion follows.
\hfill \qed
\begin{coro} \label{c95.43}
Every Fredholm sequence in $\cS_\eta (\cO_N)$ is stable.
\end{coro}
The Fredholm property of a sequence $(A_n)$ of matrices can be
expressed in terms of the singular values of the $A_n$. To
specify this remark to sequences $(A_n) \in \cS_\eta (\cO_N)$,
let
\begin{equation} \label{e95.44}
0 \le s_1(A_n) \le s_2(A_n) \le \ldots \le s_{N^n} (A_n)
= \|A_n\|
\end{equation}
be the decreasingly ordered sequence of the singular values
of $A_n$. The following is then an immediate consequence of
Corollary 6.14 in \cite{Roc10}. Note that the condition of
essential fractality, which is assumed in Corollary 6.14, is
evidently satisfied in the present context.
\begin{prop} \label{p95.45}
Let $(A_n) \in \cS_\eta (\cO_N)$. If the sequence $(A_n)$ is
stable, then the sequence $(s_1(A_n))$ is bounded below from
zero by a positive constant. If the sequence $(A_n)$ fails
to be stable, then
\[
\lim_{n \to \infty} s_k(A_n) = 0 \quad \mbox{for each} \quad
k \ge 1.
\]
\end{prop}
\paragraph{An open question.} In \cite{HRS2} there are
considered discretizations by the finite sections method of
several concrete algebras $\cA$ of operators on $l^2(\sZ^+)$.
In all cases, we observed that the associated 
quasicommutator ideal $\cJ(\cA)$ is a direct (or $c_0-$) sum of elementary ideals (i.e. ideals which are isomorphic to the ideal of the compact operators on a Hilbert space). Above we have seen that the ideal $\cJ_N$ is isomorphic to $\cC_N$ and, hence, to a corner of a tensor product of an AF-algebra by $K(l^2(\sZ^+))$. Is this the archetypal picture of a fractal irreducible ideal in
$\cF$?
\section{Spatial discretization of the extended Cuntz algebra}
\label{s7}
Here we consider the smallest $C^*$-subalgebra $\cO_N^{ext}$
of $L(l^2(\sZ^+))$ which contains the (concrete) Cuntz algebra
$\cO_N$ and the ideal $K(l^2(\sZ^+))$ of the compact operators.
We let $\cS_\eta(\cO_N^{ext})$ denote the smallest
$C^*$-subalgebra of the sequence algebra $\cF_\eta$ which
contains all finite sections sequences $(P_{N^n} A P_{N^n})
_{n \ge 0}$ of operators $A \in \cO_N^{ext}$ and set $\cS_N^{ext}
:= \cS_\eta(\cO_N^{ext})/\cG_\eta$. Further, let
\[
\cJ_N^{comp} := \{ (P_{N^n} K P_{N^n})_{n \ge 0} + \cG_\eta :
K \in K(l^2(\sZ^+)) \}.
\]
\begin{prop} \label{p96.1}
The set $\cJ_N^{comp}$ is a closed ideal of $\cS_N^{ext}$, and
every coset $a \in \cS_N^{ext}$ can be uniquely written as $b
+ k$ with $b \in \cS_N$ and $k \in \cJ_N^{comp}$.
\end{prop}
{\bf Proof.} It is evident from the above definitions that
$\cJ_N^{comp}$ is contained in $\cS_N^{ext}$, and it is easy
to see that $\cJ_N^{comp}$ is a closed ideal of $\cS_N^{ext}$
(compare the proof of Theorem 1.19 in \cite{HRS2}). Since
every sequence $(P_{N^n} K P_{N^n})_{n \ge 0}$ with a compact
operator $K$ is compact, one has $\cJ_N^{comp} \subset
\cK_\eta/\cG_\eta$. From Proposition \ref{p95.42} we thus
conclude that $\cS_N^{ext} \cap \cJ_N^{comp}$ is the zero
ideal, whence the assertion.
\hfill \qed \\[3mm]
We let $W$ be the mapping which associates with every sequence
in $\cS_\eta(\cO_N^{ext})$ its strong limit. Clearly, $W$ is a
$^*$-homomorphism which has the ideal $\cG_\eta$ in its kernel.
Thus, the quotient homomorphism
\[
\cS_N^{ext} \to L(l^2(\sZ^+)), \quad (A_n)_{n \ge 0} + \cG_\eta
\mapsto W((A_n)_{n \ge 0})
\]
is correctly defined; we denote it by $W$ again. Further we
extend the definition of the homomorphism $\widetilde{W}$ from
$\cS_N$ to the extended algebra $\cS_N^{ext}$ by writing every
element $a$ of $\cS_N^{ext}$ as $b+k$ with $b$ and $k$ as in
Proposition \ref{p96.1} and setting $\widetilde{W} (a) :=
\widetilde{W} (b)$. It is easy to check that $\widetilde{W}$ is
a $^*$-homomorphism on the extended algebra $\cS_N^{ext}$. Note
that $\cJ_N$ is in the kernel of $W$ and $\cJ_N^{comp}$ is in
the kernel of (the extended) $\widetilde{W}$.
\begin{remark}
One can also give an explicit definition of the extended
homomorphism $\widetilde{W}$. For, recall that $\widetilde{W}$
is defined on $\cS_N$ via the homomorphism $W_{00}$ which
on its hand is defined on cosets $(P_{N^n} A P_{N^n}) + \cG_N$
with $A \in \cO_N^{par}$  by the strong limit
\[
\mbox{s-lim}_{n \to \infty} R_{N^n} P_{N^n} A P_{N^n} R_{N^n}
\]
(compare (\ref{e95.13})). It is easy to check that this limit
exists as well if $A$ is a compact operator, in which case this
limit is zero. \hfill \qed
\end{remark}
\begin{theo} \label{t96.2}
A sequence $\bA = (A_n) \in \cS_\eta(\cO_N^{ext})$ is stable if
and only if the operators $W(\bA + \cG_\eta)$ and $\widetilde{W}
(\bA + \cG_\eta)$ are invertible.
\end{theo}
{\bf Proof.} It is evident that the stability of $\bA$ implies
the invertibility of $W(\bA + \cG_\eta)$ and $\widetilde{W}
(\bA + \cG_\eta)$. Conversely, assume that these operators
are invertible for a sequence $\bA \in \cS_\eta(\cO_N^{ext})$.
Set $a := \bA + \cG_\eta$ and write $a$ as $b+k$ with $b$ and
$k$ as in Proposition \ref{p96.1}. Then $\widetilde{W} (b) =
\widetilde{W} (a)$ is invertible, whence the invertibility of
$b$ by Theorem \ref{t95.38}. Thus, the sequence $\bA$ is the
sum of a stable sequence and of a sequence of the form
$(P_{N^n} K P_{N^n})$ with a compact operator $K$. Now
the assertion follows immediately from a general result on
compact perturbations; see Corollary 1.22 in \cite{HRS2}. \hfill
\qed \\[3mm]
The following is an immediate consequence of this stability
result and of Theorem 1.69 in \cite{HRS2}. Note that the
homomorphisms $W$ and $\widetilde{W}$ are fractal in the sense
of Definition 1.62 in \cite{HRS2}. This fact is evident for $W$,
and it has been checked in the proof of Corollary \ref{c95.40}
for $\widetilde{W}$.
\begin{coro} \label{c96.3}
The algebra $\cS_\eta(\cO_N^{ext})$ is fractal.
\end{coro}
As above, fractality has striking consequences for the asymptotic
behavior of some spectral quantities.
\begin{theo} \label{t96.4}
Let $(A_n)$ be a sequence in $\cS_\eta(\cO_N^{ext})$ and set $a
:= (A_n) + \cG_\eta$. Then the following set-sequences converge
with respect to the Hausdorff distance as $n \to \infty \!:$
\\[1mm]
$(a)$ $\sigma(A_n) \to \sigma (W(a)) \cup \sigma (\widetilde{W}
(a))$ if $a$ is self-adjoint;
\\[1mm]
$(b)$ $\sigma_2 (A_n) \to \sigma_2 (W(a)) \cup \sigma_2
(\widetilde{W}(a))$;
\\[1mm]
$(c)$ $\sigma^{(\varepsilon)} (A_n) \to \sigma^{(\varepsilon)}
(W(a)) \cup \sigma^{(\varepsilon)} (\widetilde{W}(a))$.
\end{theo}
Finally we consider the Fredholm properties of sequences in
$\cS_\eta(\cO_N^{ext})$. We use the notation (\ref{e95.44})
for the singular values.
\begin{theo} \label{t96.5}
Let $(A_n)$ be a sequence in $\cS_\eta(\cO_N^{ext})$ and set $a
:= (A_n) + \cG_\eta$. \\[1mm]
$(a)$ The sequence $(A_n)$ is Fredholm if and only if the
operator $\widetilde{W}(a)$ is invertible. In this case, the
operator $W(a)$ is Fredholm. \\[1mm]
$(b)$ Let $(A_n)$ be a Fredholm sequence and $k := \dim \ker
W(a)$. Then
\begin{equation} \label{e96.6}
\lim_{n \to \infty} \sigma_k (A_n) = 0 \quad \mbox{and} \quad
\liminf_{n \to \infty} \sigma_{k+1} (A_n) > 0.
\end{equation}
\end{theo}
{\bf Proof.} For assertion $(a)$, write $a$ as $b+k$ with $b$ and
$k$ as in Proposition \ref{p96.1}. If $a$ is invertible modulo
$\cJ_N^{comp}$, then $b$ is invertible modulo $\cJ_N^{comp}$;
hence $b$ is invertible due to Corollary \ref{c95.43}. Thus,
$\widetilde{W}(a) = \widetilde{W}(b)$ is invertible. Conversely,
if $\widetilde{W}(a)$ is invertible, then $\widetilde{W}(b)$
is invertible, whence the invertibility of $b$ via Theorem
\ref{t95.38}. But then $b+k = a$ is Fredholm.

Assertion $(b)$ can be derived from Theorem 6.12 in \cite{HRS2}.
But note that the assumption of a {\em standard} algebra which
is made in \cite{HRS2} is not satisfied in the present context.
A look at the proof of 6.12 in \cite{HRS2} shows that this
assumption is not really needed. A short direct proof is in \cite{Sil9}.
\hfill \qed
\section{Appendix}
The algebra $\cF$ of all bounded sequences of matrices is an algebra with polar decomposition in the following sense: A $C^*$-algebra $\cA$ has the {\em polar decomposition property} if every element $a$ of $\cA$ can be written as $a = r u$ with a unitary element $u \in \cA$ and a positive element $r \in \cA$. It is well-known from Linear Algebra that every matrix $A \in \sC^{n \times n}$ admits a polar decomposition. Thus, $\sC^{n \times n}$ owns the polar decomposition property, and so does the algebra $\cF$ and each quotient of $\cF$.

If $\cA$ is a unital $C^*$-algebra with polar decomposition property, then every element of $\cA$ which is invertible from one side is invertible from both sides. Indeed, let $b \in \cA$ be a left inverse of $a \in \cA$, and write $a$ as $a = r u$. From $ba = b r u = e$ we get $b r = u^*$ and $u b r = e$. Hence, $r$ is invertible from the left-hand side. Taking adjoints in $u b r = e$ we get the invertibility of $r$ from the right. Thus, $r$ is invertible, and so is $a$. In particular, every isometry in an
algebra with polar decomposition is unitary.
\begin{coro} \label{c97.1}
Algebras with polar decomposition cannot contain Cuntz algebras $\cO_N$ with $N \ge 2$ as subalgebras.
\end{coro}
Indeed, let $\cA$ be a $C^*$-algebra with polar decomposition, let $N \ge 2$, and assume there are isometries  $s_0, \, \ldots, \, s_{N-1} \in \cA$ with $s_0 s_0^* + \ldots + s_{N-1} s_{N-1}^* = e$. Then the $s_i$ are unitaries, whence $s_i s_i^* = e$ for every $i$. Substituting $s_i s_i^* = e$ in the Cuntz axiom gives $Ne = e$, a contradiction.
\begin{coro} \label{c97.2}
For $N \ge 2$, there are no proper closed ideal $\cJ$ of $\cF$
and no sequences $(S_n^{(0)})_{n \ge 1}, \, \ldots, \,
(S_n^{(N-1)})_{n \ge 1}$ in $\cF$ such that
\begin{equation} \label{e97.2}
\left( (S_n^{(i)})^* S_n^{(i)} - I_n \right) \in \cJ \quad
\mbox{for all} \; \; i
\end{equation}
and
\begin{equation} \label{e97.3}
\left( S_n^{(0)} (S_n^{(0)})^* + \ldots + S_n^{(N-1)} (S_n^{(N-1)})^* - I_n \right) \in \cJ.
\end{equation}
\end{coro}
This argument does not work if we replace $\cF$ by a subalgebra,
$\cB$ say, since the factors of the polar decomposition need not
belong to $\cB$. 

Note that, in contrast to Cuntz algebras, the Toeplitz algebra $\cT(C)$, which is the universal algebra generated by one isometry, can be embedded in a quotient of $\cF$. Indeed, let $\cS(\cT(C))$ stand for smallest closed subalgebra of $\cF$ which contains all sequences $(P_n A P_n)$ of finite sections of operators $A \in \cT(C)$. One can show (see Theorem 1.53 in \cite{HRS2} for instance) that every sequence $(A_n)$ in $\cS(\cT(C))$ can be uniquely written as
\[
(A_n) = (P_n T(f) P_n) + (P_n K P_n) + (R_n L R_n) + (G_n)
\]
where $T(f)$ is a Toeplitz operator with continuous generating function $f$, $K$ and $L$ are compact operators, and $(G_n)$ is a sequence tending to zero in the norm (recall the definition of the operators $R_n$ at the beginning of Section \ref{ss62}). It is not hard to check that the set
\[
\cJ := \{ (R_n L R_n) + (G_n) : K \; \mbox{compact}, \; \|G_n\| \to 0 \}
\]
forms a closed ideal of $\cS(\cT(C))$ and that the quotient $\cS(\cT(C))/\cJ$ is isomorphic to $\cT(C)$. In particular, $\cS(\cT(C))/\cJ$ contains (and is generated by) the non-unitary isometry $(P_n V P_n) + \cJ$ where $V$ is the forward shift operator on $l^2(\sZ^+)$. Note that $\cJ$ is contained in the ideal $\cK$ of the compact sequences of $\cF$ defined in Section \ref{ss33}.
{\small Author's address: \\[3mm]
Steffen Roch, Technische Universit\"at Darmstadt, Fachbereich
Mathematik, Schlossgartenstrasse 7, 64289 Darmstadt,
Germany. \\
E-mail: roch@mathematik.tu-darmstadt.de}
\end{document}